\documentclass[aos,MSNbibl,dvips]{arximspdf}

\doi{10.1214/14-AOS1235} 
\volume{42}
\issue{5}
\pubyear{2014}
\firstpage{1787}
\lastpage{1818}
\docsubty{FLA}

\makeatletter
\renewcommand{\epsilon}{\varepsilon}
\newcommand{\rrvert}{\vert}
\newcommand{\llvert}{\vert}
\newtheorem{theorem}{Theorem}[section]
\newtheorem{lemma}{Lemma}[section]
\newtheorem{corollary}{Corollary}[section]
\newproclaim{definition}{Definition}[section]

\newproclaim{comment}{Comment}[section]
\newproclaim{condition}{Condition} 
\newproclaim{conditionlow}{Condition} 
\newproclaim{conditionvc}{Condition} 
\newtheorem{algorithm}{Algorithm}
\newcommand{\Ep}{{\mathrm{E}}}
\renewcommand{\Pr}{{\mathrm{P}}}
\newcommand{\var}{\operatorname{Var}}

\newcommand{\widebar}{\overline}

\makeatother

\begin{document}
\begin{frontmatter}

\title{Anti-concentration and honest, adaptive confidence~bands}
\runtitle{Anti-concentration and confidence bands}

\begin{aug}
\author[a]{\fnms{Victor}~\snm{Chernozhukov}\thanksref{t1}\ead[label=e1]{vchern@mit.edu}},
\author[b]{\fnms{Denis}~\snm{Chetverikov}\corref{}\thanksref{t1}\ead[label=e2]{chetverikov@econ.ucla.edu}}\break
\and
\author[c]{\fnms{Kengo}~\snm{Kato}\thanksref{t2}\ead[label=e3]{kkato@e.u-tokyo.ac.jp}}
\runauthor{V. Chernozhukov, D. Chetverikov and K. Kato}
\affiliation{Massachusetts Institute of Technology,
University of California,\break  Los Angeles and University of Tokyo}
\address[a]{V. Chernozhukov\\
Department of Economics and\\
\quad Operations Research Center\\
Massachusetts Institute of Technology\\
50 Memorial Drive\\
Cambridge, Massachusetts 02142\\
USA\\
\printead{e1}}
\address[b]{D. Chetverikov\\
Department of Economics\\
University of California, Los Angeles\\
Bunche Hall, 8283\\
315 Portola Plaza\\
Los Angeles, California 90095\\
USA\\
\printead{e2}}
\address[c]{K. Kato\\
Graduate School of Economics\\
University of Tokyo\\
7-3-1 Hongo, Bunkyo-ku\\
Tokyo 113-0033\\
Japan\\
\printead{e3}}
\end{aug}
\thankstext{t1}{Supported by a National Science Foundation grant.}
\thankstext{t2}{Supported by the Grant-in-Aid for Young Scientists (B)
(25780152), the Japan Society for the Promotion of Science.}

\received{\smonth{12} \syear{2013}}
\revised{\smonth{4} \syear{2014}}

%
\begin{abstract}
Modern construction of uniform confidence bands for nonparametric
densities (and other functions) often relies on the classical
Smirnov--Bickel--Rosenblatt (SBR) condition; see, for example, Gin\'{e}
and Nickl [\textit{Probab. \mbox{Theory} Related Fields} \textbf{143} (2009) 569--596].
This condition requires the existence of a limit distribution of an
extreme value type for the supremum of a studentized empirical process
(equivalently, for the supremum of a Gaussian process with the same
covariance function as that of the studentized empirical process).
The principal contribution of this paper is to remove the need for this
classical condition. We show that a considerably weaker sufficient
condition is derived from an anti-concentration property of the
supremum of the approximating Gaussian process, and we derive an
inequality leading to such a property for separable Gaussian processes.
We refer to the new condition as a generalized SBR condition. Our new
result shows that the supremum does not concentrate too fast around any value.

We then apply this result to derive a Gaussian multiplier bootstrap
procedure for constructing honest confidence bands for nonparametric
density estimators (this result can be applied in other nonparametric
problems as well). An essential advantage of our approach is that it
applies generically even in those cases where the limit distribution of
the supremum of the studentized empirical process does not exist (or is
unknown). This is of particular importance in problems where resolution
levels or other tuning parameters have been chosen in a data-driven
fashion, which is needed for adaptive constructions of the confidence
bands. 
Finally, of independent interest is our introduction of a new,
practical version of Lepski's method, which computes the optimal,
nonconservative resolution levels via a Gaussian multiplier bootstrap method.
\end{abstract}

%
\end{frontmatter}

\section{Introduction}\label{sec1}

Let $X_{1},\ldots,X_{n}$ be i.i.d. random vectors with common unknown
density $f$ on $\mathbb{R}^d$.
We are interested in constructing confidence bands for $f$ on a subset
$\mathcal{X} \subset\mathbb{R}^d$ that are \emph{honest} to a given
class $\mathcal{F}$ of densities on $\mathbb{R}^d$. Typically,
$\mathcal{X}$ is a compact set on which $f$ is bounded away from zero,
and $\mathcal{F}$ is a class of smooth densities such as a subset of a
H\"{o}lder ball. A confidence band $\mathcal{C}_{n} = \mathcal
{C}_{n}(X_{1},\ldots,X_{n})$ is a family of random intervals
\[
\mathcal{C}_{n}:= \bigl\{ \mathcal{C}_{n}(x) =
\bigl[c_{L}(x), c_{U}(x)\bigr]\dvtx  x \in\mathcal{X} \bigr\}
\]
that contains the graph of $f$ on $\mathcal{X}$ with a guaranteed
probability. Following \cite{L89}, a band $\mathcal{C}_{n}$ is said
to be \emph{asymptotically honest with level $\alpha\in(0,1)$ for the
class $\mathcal{F}$} if
\[
\liminf_{n \to\infty} \inf_{f \in\mathcal{F}} \Pr_{f}
\bigl( f(x) \in\mathcal{C}_{n}(x),\ \forall x \in\mathcal{X} \bigr) \geq
1-\alpha.
\]
Also, we say that a band $\mathcal{C}_n$ is \textit{asymptotically
honest at a polynomial rate with level $\alpha\in(0,1)$ for the class
$\mathcal{F}$} if
%
%
\begin{equation}
\label{eqhonestatpolynomialrate} \inf_{f \in\mathcal{F}} \Pr_{f} \bigl( f(x) \in
\mathcal {C}_{n}(x),\ \forall x \in\mathcal{X} \bigr) \geq1-\alpha-C
n^{-c}
\end{equation}
for some constants $c,C>0$.

Let $\hat{f}_{n}(\cdot,l)$ be a generic estimator of $f$ with a
smoothing parameter $l$, say bandwidth or resolution level, where $l$
is chosen from a candidate set $\mathcal{L}_{n}$; see \cite
{Hardle,Tsybakov09,Wasserman06} for a textbook level introduction to
the theory of density estimation. Let $\hat{l}_{n} = \hat
{l}_{n}(X_{1},\ldots,X_{n})$ be a possibly data-dependent choice of $l$
in $\mathcal{L}_{n}$.
Denote by $\sigma_{n,f}(x,l)$ the standard deviation of $\sqrt{n}\hat
{f}_{n}(x,l)$, that is, $\sigma_{n,f}(x,l):= (n\var_{f} (\hat
{f}_{n}(x,l)))^{1/2}$. Then we consider a confidence band of the form
%
%
\begin{equation}
\mathcal{C}_{n}(x) = \biggl[ \hat{f}_{n}(x,
\hat{l}_{n}) - \frac
{c(\alpha) \sigma_{n,f}(x,\hat{l}_{n})}{\sqrt{n}}, \hat {f}_{n}(x,
\hat{l}_{n}) + \frac{c(\alpha) \sigma_{n,f}(x,\hat
{l}_{n})}{\sqrt{n}} \biggr], \label{band}
\end{equation}
where $c(\alpha)$ is a (possibly data-dependent) critical value
determined to make the confidence band to have level $\alpha$.
Generally, $\sigma_{n,f}(x,l)$ is unknown and has to be replaced by an
estimator.

A crucial point in construction of confidence bands is the computation
of the critical value $c(\alpha)$.
Assuming that $\sigma_{n,f}(x,l)$ is positive on $\mathcal{X}\times
\mathcal{L}_n$, define the stochastic process
%
%
\begin{equation}
Z_{n,f}(v):= Z_{n,f}(x,l):= \frac{\sqrt{n}(\hat{f}_{n}(x,l)-\Ep
_{f}[\hat{f}_{n}(x,l)])}{\sigma_{n,f}(x,l)}, \label{student}
\end{equation}
where\vspace*{1pt} $v=(x,l)\in\mathcal{X}\times\mathcal{L}_n=:\mathcal{V}_n$.
We refer to $Z_{n,f}$ as a ``studentized process.''
If, for the sake of simplicity, the bias $|f(x)-\Ep_{f}[ \hat
{f}_{n}(x,l) ]_{l=\hat{l}_{n}}|$ is sufficiently small compared to
$\sigma_{n,f}(x,\hat{l}_{n})$, then
\begin{eqnarray*}
\Pr_{f} \bigl(f(x) \in\mathcal{C}_{n}(x),\ \forall x \in
\mathcal {X} \bigr) &\approx& \Pr_{f} \Bigl( \sup_{x \in\mathcal{X}}
\bigl\llvert Z_{n,f}(x,\hat{l}_{n}) \bigr\rrvert \leq c(
\alpha) \Bigr)
\\
&\geq& \Pr_{f} \Bigl( \sup_{v \in\mathcal{V}_{n}} \bigl\llvert
Z_{n,f}(v) \bigr\rrvert \leq c(\alpha) \Bigr),
\end{eqnarray*}
so that band (\ref{band}) will be of level $\alpha\in(0,1)$ by taking
%
%
\begin{equation}
c(\alpha) = (1-\alpha)\mbox{-quantile of } \| Z_{n,f} \|_{\mathcal
{V}_{n}}:= \sup_{v \in\mathcal{V}_{n}}\bigl| Z_{n,f}(v)\bigr|.
\end{equation}
The critical value $c(\alpha)$, however, is infeasible since the
finite sample distribution of the process $Z_{n,f}$ is unknown.
Instead, we estimate the $(1-\alpha)$-quantile of $\| Z_{n,f} \|
_{\mathcal{V}_{n}}$.

Suppose that one can find an appropriate centered Gaussian process
$G_{n,f}$ indexed by $\mathcal{V}_{n}$ with known or estimable
covariance structure such that
$\Vert Z_{n,f}\Vert_{\mathcal{V}_n}$ is close to $\Vert G_{n,f}\Vert
_{\mathcal{V}_n}$. Then we may approximate the $(1-\alpha)$-quantile
of $\| Z_{n,f} \|_{\mathcal{V}_{n}}$ by
\[
c_{n,f}(\alpha):= (1-\alpha)\mbox{-quantile of }\| G_{n,f} \|
_{\mathcal{V}_{n}}.
\]
Typically, one computes or approximates $c_{n,f}(\alpha)$ by one of
the following two methods:

\begin{longlist}[(2)]
\item[(1)] Analytical method: derive analytically an approximated
value of $c_{n,f}(\alpha)$, by using an explicit limit distribution or
large deviation inequalities.
\item[(2)] Simulation method: simulate the Gaussian process $G_{n,f}$
to compute $c_{n,f}(\alpha)$ numerically, by using, for example, a
multiplier method.
\end{longlist}

The main purpose of this paper is to introduce a general approach to
establishing the validity of the so-constructed confidence band.
Importantly, our analysis does not rely on the existence of an explicit
(continuous) limit distribution of any kind, which is a major
difference from the previous literature.
For the density estimation problem, if $\mathcal{L}_{n}$ is a
singleton, that is, the smoothing parameter is chosen
deterministically, the existence of such a continuous limit
distribution, which is typically a Gumbel distribution, has been
established for
convolution kernel density estimators and \emph{some} wavelet projection
kernel density estimators; see \cite
{S50,BR73,GKS04,GN10,B11a,B11b,GGM11}. We refer to the existence of the
limit distribution as the Smirnov--Bickel--Rosenblatt (SBR) condition.
However, the SBR condition has not been obtained for other density
estimators such as nonwavelet projection kernel estimators based, for
example, on Legendre polynomials or Fourier series.
In addition, to guarantee the existence of a continuous limit
distribution often requires more stringent regularity conditions than a
Gaussian approximation itself.
More importantly, if $\mathcal{L}_{n}$ is not a singleton, which is
typically the case when $\hat{l}_{n}$ is data-dependent, and so the
randomness of $\hat{l}_{n}$ has to be taken into account, it is often
hard to determine an exact limit behavior of $\| G_{n,f} \|_{\mathcal{V}_{n}}$.

We thus take a different route and significantly generalize the SBR
condition. Our key ingredient is the \emph{anti-concentration} property
of suprema of Gaussian processes that shows that suprema of Gaussian
processes do not concentrate too fast. To some extent, this is a
reverse of numerous concentration inequalities for Gaussian processes.
In studying the effect of approximation and estimation errors on the
coverage probability, it is required to know how the random variable
$\| G_{n,f} \|_{\mathcal{V}_{n}}:= \sup_{v \in\mathcal{V}_{n}} |
G_{n,f}(v) |$ concentrates or ``anti-concentrates'' around, say, its
$(1-\alpha)$-quantile.
It is not difficult to see that $\| G_{n,f} \|_{\mathcal{V}_{n}}$
itself has a continuous distribution, so that \emph{with keeping $n$
fixed}, the probability that
$\| G_{n,f} \|_{\mathcal{V}_{n}}$ falls into the interval with center
$c_{n,f}(\alpha)$ and radius $\epsilon$ goes to $0$ as $\epsilon\to0$.
However, what we need to know is the behavior of those probabilities
when $\epsilon$ depends on $n$ and $\epsilon= \epsilon_{n} \to0$.
In other words, bounding explicitly ``anti-concentration''
probabilities for suprema of Gaussian processes is desirable.
We will first establish bounds on the L\'{e}vy concentration function
(see Definition~\ref{levy}) for suprema of Gaussian processes and then
use these bounds to quantify the effect of approximation and estimation
errors on the finite sample coverage probability. We say that a \emph{generalized} SBR condition or simply an anti-concentration condition
holds if $\| G_{n,f} \|_{\mathcal{V}_{n}}$ concentrates sufficiently
slowly, so that this effect is sufficiently small to yield
asymptotically honest confidence bands.

As a substantive application of our results, we consider the problem of
constructing honest \emph{adaptive} confidence bands based on either
convolution or wavelet projection kernel density estimators in H\"
{o}lder classes $\mathcal{F}\subset\bigcup_{t\in[\underline{t},\bar
{t}]}\Sigma(t,L)$ for some $0<\underline{t}<\bar{t}<\infty$ where
$\Sigma(t,L)$ is the H\"{o}lder ball of densities with radius $L$ and
smoothness level $t$. Following \cite{CaiLow04}, we say the confidence
band $\mathcal{C}_n$ is \emph{adaptive} if for every $t,\varepsilon>0$
there exists $C>0$ such that for all $n \geq1$,
\[
\sup_{f\in\mathcal{F}\cap\Sigma(t,L)}\Pr_f \Bigl(\sup_{x\in
\mathcal{X}}
\lambda\bigl(\mathcal{C}_{n}(x)\bigr)>Cr_n(t) \Bigr)\leq
\varepsilon,
\]
where $\lambda$ denotes the Lebesgue measure on $\mathbb{R}$ and
$r_n(t):=(\log n/n)^{t/(2t+d)}$, the minimax optimal rate of
convergence for estimating a density $f$ in the function class $\Sigma
(t,L)$ in the sup-metric $d_{\infty}(\hat{f},f)=\sup_{x\in\mathcal
{X}}|\hat{f}(x)-f(x)|$. We use Lepski's method \cite
{Lepski91,Birge2001} to find an adaptive value of the smoothing
parameter. Here our contribution is to introduce a \textit{Gaussian
multiplier bootstrap implementation} of Lepski's method. This is a
practical proposal since previous implementations relied on
conservative (one-sided) maximal inequalities and are not necessarily
recommended for practice; see, for example, \cite{GN09b} for a discussion.

We should also emphasize that our techniques can also be used for
constructing honest and/or adaptive confidence bands in many other
nonparametric problems, but in this paper we focus on the density
problem for the sake of clarity. Our techniques [anti-concentration of
separable Gaussian processes (Theorem~\ref{thm2}), and coupling
inequalities (Theorems~\ref{thmcouplinginequality1}
and~\ref{thmcouplinginequality2})] are of particular importance in non-Donsker
settings since they allow us to prove validity of the Gaussian
multiplier bootstrap for approximating distributions of suprema of
sequences of empirical processes of VC type function classes where the
metric entropy of the process may increase with $n$. Thus these
techniques may be important in many nonparametric problems. 
For example, applications of our anti-concentration bounds can be found
in \cite{CLR12} and \cite{Chetverikov2012}, which consider the
problems of nonparametric inference on a minimum of a function and
nonparametric testing of qualitative hypotheses about functions, respectively.



\subsection{Related references}

Confidence bands in nonparametric estimation have been extensively
studied in the literature.
A classical approach, which goes back to \cite{S50} and \cite{BR73},
is to use explicit limit distributions of normalized suprema of
studentized processes.
A ``Smirnov--Bickel--Rosenblatt type limit theorem'' combines Gaussian
approximation techniques and extreme value theory for Gaussian processes.
It was argued that the convergence to normal extremes is rather slow
\mbox{despite} that the Gaussian approximation is relatively fast \cite{H91}.
To improve the finite sample coverage, bootstrap is often used in
construction of confidence bands; see \mbox{\cite{CV03,BDHM07}}.
However, to establish the validity of bootstrap confidence bands,
researchers relied on the existence of continuous limit distributions
of normalized suprema of original studentized processes.
In the deconvolution density estimation problem, Lounici and Nickl
\cite{LN11} considered confidence bands without using Gaussian
\mbox{approximation}. In the current density estimation problem, their idea
reads as bounding the deviation probability of $\| \hat{f}_{n} - \Ep[
\hat{f}_{n}(\cdot) ] \|_{\infty}$ by using Talagrand's \cite{T96}
inequality and replacing the expected supremum by the Rademacher
average. Such a construction is indeed general and applicable to many
other problems, but is likely to be more conservative than our construction.

\subsection{Organization of the paper}
In the next section, we give a
new anti-concentration inequality for suprema of Gaussian processes.
Section~\ref{secanalysisunderhighlevelconditions} contains a theory of
generic confidence band construction under high-level conditions. These
conditions are easily satisfied both for convolution and projection
kernel techniques under mild primitive assumptions, which are also
presented in Section~\ref{secanalysisunderhighlevelconditions}.
Section~\ref{secverifyingconditions} is devoted to constructing honest
adaptive confidence bands in H\"{o}lder classes. Finally, most proofs
are contained in the \hyperref[app]{Appendix}, and some proofs and discussions are put
into the supplemental material \cite{CCK2014}.

\subsection{Notation}
In what follows, constants $c,C,c_{1},C_{1},c_{2},C_{2},\ldots$ are
understood to be positive and independent of $n$. The values of $c$ and
$C$ may change at each appearance but constants
$c_{1},C_{1},c_{2},C_{2},\ldots$ are fixed.
Throughout the paper, $\mathbb{E}_n[\cdot]$ denotes the
average over index $1 \leq i \leq n$, that is, it simply abbreviates
the notation $n^{-1} \sum_{i=1}^n[\cdot]$. For example, $\mathbb
{E}_n[ g(X_{i}) ] = n^{-1} \sum_{i=1}^{n} g(X_{i})$.
For a set $T$, denote by $\ell^{\infty}(T)$ the set of all bounded
functions, that is, all functions $z\dvtx  T \to\mathbb{R}$ such that
\[
\| z\|_{T}:= \sup_{t \in T} \bigl| z(t) \bigr| < \infty.
\]
Moreover, for a generic function $g$, we also use the notation $\| g \|
_{\infty}:=\break  \sup_{x} | g(x) |$ where the supremum is taken over the
domain of $g$.
For two random variables $\xi$ and $\eta$, we write $\xi\stackrel
{d}{=} \eta$ if they share the same distribution.
The standard Euclidean norm is denoted by $| \cdot|$.

\section{Anti-concentration of suprema of Gaussian processes}
\label{secanti-concentration}
The main purpose of this section is to derive an upper bound on the
\emph{L\'{e}vy concentration function} for suprema of separable Gaussian
processes, where the terminology is adapted from \cite{RV09}. Let
$(\Omega,\mathcal{A},\Pr)$ be the underlying (complete) probability space.
%
\begin{definition}[(L\'{e}vy concentration function)]\label{levy}
Let $Y=(Y_t)_{t \in T}$ be a separable stochastic process indexed by a
semimetric space $T$. For all $x\in\mathbb{R}$ and $\epsilon\geq0$, let
%
%
\begin{equation}
p_{x,\epsilon} (Y):= \Pr \Bigl(\Bigl\llvert \sup_{t\in T}Y_t-x
\Bigr\rrvert \leq \epsilon \Bigr) \label{anti}.
\end{equation}
Then the \emph{L\'{e}vy concentration function} of $\sup_{t \in T} Y_t$
is defined for all $\epsilon\geq0$ as
%
%
\begin{equation}
p_{\epsilon}(Y):= \sup_{x \in\mathbb{R}} p_{x,\epsilon}(Y).
\label{concen}
\end{equation}
Likewise, define $p_{x,\epsilon}(|Y|)$ by (\ref{anti}) with $\sup_{t
\in T} Y_t$ replaced by $\sup_{t \in T} | Y_t |$, and define
$p_{\epsilon}(|Y|)$ by (\ref{concen}) with $p_{x,\epsilon}(Y)$
replaced by $p_{x,\epsilon}(|Y|)$.
\end{definition}

Let $X=(X_{t})_{t \in T}$ be a separable Gaussian process indexed by a
semimetric space $T$ such that $\Ep[ X_{t} ] = 0$ and $\Ep[ X^{2}_{t}
] = 1$ for all $t \in T$.
Assume that $\sup_{t \in T}X_{t} < \infty$ a.s.
Our aim here is to obtain a qualitative bound on the concentration
function $p_{\epsilon}(X)$. In a trivial example where $T$ is a
singleton, that is, $X$ is a real standard normal random variable, it
is immediate to see that $p_{\epsilon}(X) \asymp\epsilon$ as
$\epsilon\to0$.
A nontrivial case is that when $T$ is not a singleton, and both $T$
and $X$ are indexed by $n = 1,2,\ldots,$ that is, $T=T_{n}$ and
$X=X^{n}=(X_{n,t})_{t \in T_{n}}$, and the complexity of the set $\{
X_{n,t}\dvtx  t \in T_{n} \}$ [in $L^{2}(\Omega,\mathcal{A},\Pr)$] is
increasing in $n$. In such a case, it is typically not known whether
$\sup_{t \in T_{n}} X_{n,t}$ has a limiting distribution as $n \to
\infty$, and therefore it is not trivial at all whether, for any
sequence $\epsilon_{n} \to0$, $p_{\epsilon_{n}} (X^{n}) \to0$ as $n
\to\infty$.

The following is the first main result of this paper.

%
\begin{theorem}[(Anti-concentration for suprema of separable Gaussian processes)]\label{thm2}
Let $X=(X_{t})_{t \in T}$ be a separable Gaussian process indexed by a
semimetric space $T$ such that $\Ep[ X_{t} ] = 0$ and $\Ep[ X^{2}_{t}
] = 1$ for all $t \in T$. Assume that $\sup_{t \in T} X_{t} < \infty$
a.s. Then $a(X):= \Ep[ \sup_{t \in T} X_{t}] \in[0,\infty)$ and
%
%
\begin{equation}
p_{\epsilon} (X) \leq4 \epsilon \bigl(a(X) + 1 \bigr), \label
{eqanti-concentrationmain}
\end{equation}
for all $\epsilon\geq0$.
\end{theorem}

The similar conclusion holds for the concentration function of $\sup_{t \in T} |X_{t}|$.

\begin{corollary}\label{cor1}
Let $X=(X_{t})_{t \in T}$ be a separable Gaussian process indexed by a
semimetric space $T$ such that $\Ep[ X_{t} ] = 0$ and $\Ep[ X^{2}_{t}
] = 1$ for all $t \in T$. Assume that $\sup_{t \in T} X_{t} < \infty$
a.s. Then $a(|X|):= \Ep[ \sup_{t \in T} |X_{t}|] \in[\sqrt{2/\pi
},\infty)$ and
%
%
\begin{equation}
p_{\epsilon} \bigl(|X|\bigr) \leq4\epsilon\bigl(a\bigl(|X|\bigr) + 1\bigr), \label
{eqanti-concentrationmaincor}
\end{equation}
for all $\epsilon\geq0$.
\end{corollary}
We refer to (\ref{eqanti-concentrationmain}) and (\ref
{eqanti-concentrationmaincor}) as anti-concentration inequalities
because they show that suprema of separable Gaussian processes can not
concentrate too fast. The proof of Theorem~\ref{thm2} and Corollary
\ref{cor1} follows by extending the results in \cite{CCK2012b} where
we derived anti-concentration inequalities for maxima of Gaussian
random vectors. See the \hyperref[app]{Appendix} for a detailed exposition.

\section{Generic construction of honest confidence bands}\label{secanalysisunderhighlevelconditions}

We go back to the analysis of confidence bands.
Recall that we consider the following setting. We observe i.i.d. random
vectors $X_{1},\ldots,X_{n}$ with common unknown density $f \in
\mathcal{F}$ on $\mathbb{R}^d$, where $\mathcal{F}$ is a nonempty
subset of densities on $\mathbb{R}^d$. We denote by $\Pr_f$ the
probability distribution corresponding to the density $f$.
We first state the result on the construction of honest confidence
bands under certain high-level conditions and then show that these
conditions hold for most commonly used kernel density estimators.

\subsection{Main result}
Let $\mathcal{X}\subset\mathbb{R}^d$ be a set of interest.
Let $\hat{f}_{n}(\cdot,l)$ be a generic estimator of $f$ with a
smoothing parameter $l \in\mathcal{L}_{n}$ where $\mathcal{L}_n$ is
the candidate set.
Denote by $\sigma_{n,f}(x,l)$ the standard deviation of $\sqrt{n}\hat
{f}_{n}(x,l)$. We assume that $\sigma_{n,f}(x,l)$ is positive on
$\mathcal{V}_{n}:= \mathcal{X} \times\mathcal{L}_{n}$ for all $f
\in\mathcal{F}$.
Define the studentized process $Z_{n,f} = \{ Z_{n,f}(v)\dvtx  v=(x,l) \in
\mathcal{V}_{n} \}$ by (\ref{student}).
Let
\[
W_{n,f}:=\|Z_{n,f}\|_{\mathcal{V}_n}
\]
denote the supremum of the studentized process. We assume that
$W_{n,f}$ is a well-defined random variable.
Let $c_{1},C_{1}$ be some positive constants. 
We will assume the following high-level conditions.

\renewcommand{\thecondition}{H\arabic{condition}}
%
\begin{condition}[(Gaussian approximation)]\label{con1}
For every $f \in\mathcal{F}$, there exists (on a possibly enriched
probability space) a sequence of random variables $W^0_{n,f}$ such that
\textup{(i)}~$W^0_{n,f}\stackrel{d}{=}\|G_{n,f}\|_{\mathcal{V}_n}$ where
$G_{n,f}=\{G_{n,f}(v)\dvtx v\in\mathcal{V}_n\}$ is a tight Gaussian
random element in $\ell^{\infty}(\mathcal{V}_n)$ with $\Ep[
G_{n,f}(v) ] = 0, \Ep[G_{n,f}(v)^2]=1$ for all $v\in\mathcal{V}_n$,
and $\Ep[\Vert G_{n,f}\Vert_{\mathcal{V}_n}] \leq C_1\sqrt{\log n}$;
and moreover \textup{(ii)}
%
%
\begin{equation}
\label{eqcon1} \sup_{f \in\mathcal{F}} \Pr_f\bigl(
\bigl|W_{n,f} - W^0_{n,f}\bigr| > \epsilon _{1n}
\bigr) \leq\delta_{1n},
\end{equation}
where $\epsilon_{1n}$ and $\delta_{1n}$ are some sequences of
positive numbers bounded from above by~$C_1n^{-c_1}$.
\end{condition}

Analysis of uniform confidence bands often relies on the classical
Smirnov--Bickel--Rosenblatt (SBR) condition that states that for some
sequences $A_n$ and~$B_n$,
%
%
\begin{equation}
\label{eqSBRcondition} A_n\bigl(\Vert G_{n,f}\Vert_{\mathcal{V}_n}-B_n\bigr)
\stackrel{d} {\rightarrow }Z,\qquad\mbox{as } n\rightarrow\infty,
\end{equation}
where $Z$ is a Gumbel random variable; see, for example, \cite{GN10}.
Here both $A_n$~and~$B_n$ are typically of order $\sqrt{\log n}$.
However, this condition is often difficult to verify. Therefore, we
propose to use a weaker condition (recall the definition of the L\'
{e}vy concentration function given in Definition~\ref{levy}):

%
\begin{condition}[(Anti-concentration or generalized SBR condition)]\label{con2}
For any sequence $\epsilon_n$ of positive numbers, we have
\begin{eqnarray*}
&&\mbox{\textup{(a)}\quad} \sup_{f \in\mathcal{F}} p_{\epsilon_n}\bigl(|G_{n,f}|\bigr)
\rightarrow0\qquad\mbox{if }\epsilon_n\sqrt{\log n}\rightarrow 0\quad
\mbox{or}
\\
&&\mbox{\textup{(b)}\quad} \sup_{f \in\mathcal{F}} p_{\epsilon_n}\bigl(|G_{n,f}|\bigr)
\leq C_1\epsilon_n\sqrt{\log n}.
\end{eqnarray*}
\end{condition}
Note that Condition~\ref{con2}(a) follows trivially from Condition~\ref{con2}(b). In turn, under Condition~\ref{con1}, Condition~\ref{con2}(b) is a
simple consequence of Corollary~\ref{cor1}.
Condition~\ref{con2}(a) (along with Conditions~\ref{con1} and~\ref{con3}--\ref{con6} below) is sufficient to show
that the confidence bands are asymptotically honest,
but we will use Condition~\ref{con2}(b) to show that the confidence
bands are asymptotically honest at a polynomial rate.
We refer to Condition~\ref{con2} as a generalized SBR condition
because Condition~\ref{con2}(a) holds if (\ref{eqSBRcondition}) holds with
$A_n$ of order $\sqrt{\log n}$. An advantage of Condition~\ref{con2}
in comparison with the classical condition (\ref{eqSBRcondition}) is
that Condition~\ref{con2} follows easily from Corollary~\ref{cor1}.

Let $\alpha\in(0,1)$ be a fixed constant (confidence level).
Recall that $c_{n,f}(\alpha)$ is the $(1-\alpha)$-quantile of the
random variable $\| G_{n,f} \|_{\mathcal{V}_{n}}$.
If $G_{n,f}$ is pivotal, that is, independent of $f$, $c_{n,f}(\alpha)
= c_{n}(\alpha)$ can be directly computed, at least numerically.
Otherwise, we have to approximate or estimate $c_{n,f}(\alpha)$.
Let $\hat{c}_{n}(\alpha)$ be an estimator or approximated value of
$c_{n,f}(\alpha)$, where we assume that $\hat{c}_{n}(\alpha)$ is
nonnegative [which is reasonable since $c_{n,f}(\alpha)$ is
nonnegative]. The following is concerned with a generic regularity
condition on the accuracy of the estimator $\hat{c}_{n}(\alpha)$.%

%
\begin{condition}[{[Estimation error of $\hat{c}_{n}(\alpha)$]}]\label{con3}
For some sequences $\tau_{n}$, $\epsilon_{2n}$, and $\delta_{2n}$ of
positive numbers bounded from above by $C_1n^{-c_1}$, we have
\begin{eqnarray*}
&&\mbox{\textup{(a)}\quad} \sup_{f \in\mathcal{F}} \Pr_{f} \bigl( \hat
{c}_{n}(\alpha) < c_{n,f}(\alpha+ \tau_{n}) -
\epsilon_{2n} \bigr) \leq\delta_{2n}\quad\mbox{and}
\\
&&\mbox{\textup{(b)}\quad} \sup_{f \in\mathcal{F}} \Pr_{f} \bigl( \hat
{c}_{n}(\alpha) > c_{n,f}(\alpha- \tau_{n}) +
\epsilon_{2n} \bigr) \leq\delta_{2n}.
\end{eqnarray*}
\end{condition}

In the next subsection, we shall verify this condition for the
estimator $\hat{c}_{n}(\alpha)$ based upon the Gaussian multiplier
bootstrap method. Importantly, in this condition, we introduce the
sequence $\tau_n$ and compare $\hat{c}_{n}(\alpha)$ with
$c_{n,f}(\alpha+ \tau_{n})$ and $c_{n,f}(\alpha- \tau_{n})$ instead
of directly comparing it with $c_{n,f}(\alpha)$, which considerably
simplifies verification of this condition. With $\tau_n=0$ for all
$n$, we would need to have an upper bound on $c_{n,f}(\alpha
)-c_{n,f}(\alpha+\tau_n)$ and $c_{n,f}(\alpha-\tau
_n)-c_{n,f}(\alpha)$, which might be difficult to obtain in general.


The discussion in the \hyperref[sec1]{Introduction} presumes that $\sigma_{n,f}(x,l)$
were known, but of course it has to be replaced by a suitable estimator
in practice.
Let $\hat{\sigma}_{n}(x,l)$ be a generic estimator of $\sigma
_{n,f}(x,l)$. Without loss of generality, we may assume that $\hat
{\sigma}_{n}(x,l)$ is nonnegative.
Condition~\ref{con4} below states a high-level assumption on the
estimation error of $\hat{\sigma}_{n}(x,l)$. Verifying Condition
\ref{con4} is rather standard for specific examples.

%
\begin{condition}[{[Estimation error of $\hat{\sigma}_{n}(\cdot)$]}]\label{con4}
For some sequences $\epsilon_{3n}$ and $\delta_{3n}$ of positive
numbers bounded from above by $C_1n^{-c_1}$,
\[
\sup_{f \in\mathcal{F}} \Pr_{f} \biggl( \sup
_{v \in\mathcal
{V}_{n}} \biggl\llvert \frac{\hat{\sigma}_{n} (v)}{\sigma_{n,f}(v)} - 1 \biggr\rrvert >
\epsilon_{3n} \biggr) \leq\delta_{3n}.
\]
\end{condition}

We now consider strategies to deal with the bias term. We consider two
possibilities.
The first possibility is to control the bias explicitly, so that the
confidence band contains the bias controlling term. This construction
is inspired by \cite{B11a}. The advantage of this construction is that
it yields the confidence band the length of which shrinks at the
minimax optimal rate with no additional inflating terms; see Theorem
\ref{propverifyingconditionH4} below. The disadvantage, however, is
that this construction yields a conservative confidence band in terms
of coverage probability. We consider this strategy in Conditions~\ref
{con5} and~\ref{con6} and Theorem~\ref{prop1}.
The other possibility is to undersmooth, so that the bias is
asymptotically negligible, and hence the resulting confidence band
contains no bias controlling terms. This is an often used strategy;
see, for example, \cite{GN10}. The advantage of this construction is
that it sometimes yields an exact (nonconservative) confidence band,
so that the confidence band covers the true function with probability
$1-\alpha$ asymptotically exactly; see Corollary~\ref
{corundersmoothing} below. The disadvantages, however, are that this
method yields the confidence band that shrinks at the rate slightly
slower than the minimax optimal rate, and that is centered around a
nonoptimal estimator. We consider the possibility of undersmoothing in
Corollary~\ref{corundersmoothing} below. Note that Conditions~\ref{con5} and~\ref{con6} below are not assumed in Corollary~\ref
{corundersmoothing}.

We now consider the first possibility, that is, we assume that the
smoothing parameter $\hat{l}_n:=\hat{l}_n(X_1,\ldots,X_n)$, which is
allowed to depend on the data, is chosen so that the bias can be
controlled sufficiently well. Specifically, for all $l\in\mathcal
{L}_n$, define
\[
\Delta_{n,f}(l):= \sup_{x \in\mathcal{X}} \frac{\sqrt{n}| f(x) -
\Ep_{f}[\hat{f}_{n}(x,l)] |}{\sigma_{n}(x,l)}.
\]
We assume that there exists a sequence of random variables $c_n'$,
which are known or can be calculated via simulations, that control
$\Delta_{n,f}(\hat{l}_n)$. In particular, the theory in the next
subsection assumes that $c_n'$ is chosen as a multiple of the estimated
high quantile of the supremum of certain Gaussian process.

%
\begin{condition}[{[Bound on $\Delta_{n,f}(\hat{l}_n)$]}]\label{con5}
For some sequence $\delta_{4n}$ of positive numbers bounded from above
by $C_1n^{-c_1}$,
\[
\sup_{f \in\mathcal{F}}\Pr_{f} \bigl( \Delta_{n,f}(
\hat{l}_{n}) > c_n' \bigr) \leq
\delta_{4n}.
\]
\end{condition}

In turn, we assume that $c_n'$ can be controlled by $u_n\sqrt{\log n}$
where $u_n$ is a sequence of nonnegative positive numbers. Typically,
$u_n$ is either a bounded or slowly growing sequence; see, for example,
our construction under primitive conditions in the next section.
%
\begin{condition}[(Bound on $c_n'$)]\label{con6}
For some sequences $\delta_{5n}$ and $u_n$ of positive numbers where
$\delta_{5n}$ is bounded from above by $C_1n^{-c_1}$,
\[
\sup_{f\in\mathcal{F}}\Pr_f \bigl(c_n'>u_n
\sqrt{\log n} \bigr)\leq\delta_{5n}.
\]
\end{condition}
When $\mathcal{L}_n$ is a singleton, conditions like Conditions~\ref{con5} and
\ref{con6} have to be assumed. When $\mathcal{L}_n$ contains more
than one element, that is, we seek for an adaptive procedure,
verification of Conditions~\ref{con5} and~\ref{con6} is
nontrivial. In Section~\ref{secverifyingconditions}, we provide an
example of such analysis.

We consider the confidence band $\mathcal{C}_{n} = \{ \mathcal
{C}_{n}(x)\dvtx  x \in\mathcal{X} \}$ defined by
%
%
\begin{equation}
\mathcal{C}_{n}(x):= \bigl[\hat{f}_{n}(x,
\hat{l}_{n}) - s_n(x,\hat{l}_n),
\hat{f}_{n}(x,\hat{l}_{n}) + s_n(x,
\hat{l}_n) \bigr], \label{cband}
\end{equation}
where
%
%
\begin{equation}
\label{eqsn} s_n(x,\hat{l}_n):=\bigl(
\hat{c}_{n}(\alpha)+c_n'\bigr)\hat{\sigma
}_{n}(x,\hat{l}_{n})/\sqrt{n}.
\end{equation}

Define
\begin{eqnarray*}
\bar{\epsilon}_{n,f} &:=& \epsilon_{1n} +
\epsilon_{2n} + \epsilon _{3n}\bigl(c_{n,f}(
\alpha)+u_n\sqrt{\log n}\bigr),
\\
\delta_{n}&:=& \delta_{1n} + \delta_{2n}+
\delta_{3n}+\delta _{4n}+\delta_{5n}.
\end{eqnarray*}

We are now in position to state the main result of this section. Recall
the definition of L\'{e}vy concentration function (Definition~\ref{levy}).
%
\begin{theorem}[(Honest generic confidence bands)]\label{prop1}
Suppose that Conditions~\ref{con1} and~\ref{con3}--\ref{con6} are
satisfied.
Then
%
%
\begin{equation}
\label{eqresulthonestnonasymptotic} \inf_{f\in\mathcal{F}}\Pr_{f}( f \in
\mathcal{C}_{n} ) \geq(1- \alpha) - \delta_{n} -
\tau_{n} - p_{\bar{\epsilon}_{n,f}}\bigl(|G_{n,f}|\bigr).
\end{equation}
If, in addition, Condition~\ref{con2}\textup{(a)} is satisfied and $\epsilon
_{3n}u_n\sqrt{\log n} \leq C_{1} n^{-c_{1}}$, then
%
%
\begin{equation}
\label{eqresultweak} \liminf_{n\rightarrow\infty}\inf_{f\in\mathcal{F}}
\Pr_f(f\in \mathcal{C}_{n})\geq1-\alpha,
\end{equation}
and if, in addition, Condition~\ref{con2}\textup{(b)} is satisfied, then
%
%
\begin{equation}
\label{eqresultstrong} \inf_{f\in\mathcal{F}}\Pr_f(f\in
\mathcal{C}_{n})\geq1-\alpha-Cn^{-c},
\end{equation}
where $c$ and $C$ are constants depending only on $\alpha$, $c_1$ and $C_1$.
\end{theorem}

The confidence band defined in (\ref{cband}) and (\ref{eqsn}) is
constructed so that the bias $\Delta_{n,f}(\hat{l}_n)$ is
controlled\vadjust{\goodbreak}
explicitly via the random variable $c_n'$. Alternatively, one can
choose to undersmooth so that the bias is negligible asymptotically. To
cover this possibility, we note that it follows from the proof of
Theorem~\ref{prop1} that if \mbox{$u_n\log n\to0$} or $u_n\log n\leq
C_1n^{-c_1}$, then conclusions (\ref{eqresultweak}) or (\ref
{eqresultstrong}) of Theorem~\ref{prop1} continue to hold,
respectively, with $s_n(x,\hat{l}_n)$ in (\ref{eqsn}) replaced by
$\hat{c}_n(\alpha)\hat{\sigma}_n(x,\hat{l}_n)/\sqrt{n}$. Thus,
obtaining the asymptotically honest at a polynomial rate confidence
band requires polynomial undersmoothing ($u_n\log n\leq C_1n^{-c_1}$),
but on the other hand, logarithmic undersmoothing ($u_n\log n\to0$)
suffices if polynomial rate is not required. Moreover, if $\mathcal
{L}_n$ is a singleton, it is possible to show that the confidence band
is asymptotically exact, with a polynomial convergence rate (\ref
{eqmainresultcorexactstrong}) under the condition $u_n\log n\leq
C_1n^{-c_1}$. We collect these observations into the following
corollary, the detailed proof of which can be found in the supplemental
material \cite{CCK2014}.
%
\begin{corollary}[(Honest\vspace*{2pt} generic confidence bands with undersmoothing)]\label{corundersmoothing}
Consider the confidence band $\widetilde{\mathcal{C}}_n=\{\widetilde
{\mathcal{C}}_n(x)\dvtx x\in\mathcal{X}\}$ defined by
\[
\widetilde{\mathcal{C}}_n(x):= \bigl[\hat{f}_{n}(x,
\hat{l}_{n}) - \tilde{s}_n(x,\hat{l}_n),
\hat{f}_{n}(x,\hat{l}_{n}) + \tilde{s}_n(x,
\hat{l}_n) \bigr],
\]
where
\[
\tilde{s}_n(x,\hat{l}_n):=\hat{c}_n(
\alpha)\hat{\sigma }_n(x,\hat{l}_n)/\sqrt{n}.
\]
Suppose that Conditions~\ref{con1},~\ref{con3} and~\ref{con4} are
satisfied. In addition, assume that for some sequences $\delta_{6n}$
and $u_n$ of positive numbers,
%
%
\begin{equation}
\label{eqbiascontrolmaincor} \sup_{f\in\mathcal{F}}\Pr_f \bigl(
\Delta_{n,f}(\hat {l}_n)>u_n\sqrt{\log n} \bigr)
\leq\delta_{6n},
\end{equation}
where $\delta_{6n}$ is bounded from above by $C_1n^{-c_1}$. If
Condition~\ref{con2}\textup{(a)} holds and $u_n\log n\to0$, then
%
%
\begin{equation}
\label{eqmainresultcorweak} \liminf_{n\rightarrow\infty}\inf_{f\in\mathcal{F}}
\Pr_f(f\in \widetilde{\mathcal{C}}_{n})\geq1-\alpha.
\end{equation}
If Condition~\ref{con2}\textup{(b)} holds and $u_n\log n\leq C_1n^{-c_1}$, then
%
%
\begin{equation}
\label{eqmainresultcorstrong} \inf_{f\in\mathcal{F}}\Pr_f(f\in\widetilde{\mathcal{C}}_{n})\geq 1-\alpha-Cn^{-c}.
\end{equation}
Moreover, assume in addition that $\mathcal{L}_n$ is a singleton. If
Condition~\ref{con2}\textup{(a)} holds and $u_n\log n\to0$, then
%
%
\begin{equation}
\label{eqmainresultcorexactweak} \lim_{n\rightarrow\infty}\sup_{f\in\mathcal{F}}\bigl
\llvert \Pr_f(f\in \widetilde{\mathcal{C}}_{n})-(1-\alpha)
\bigr\rrvert = 0.
\end{equation}
If Condition~\ref{con2}\textup{(b)} and $u_n\log n\leq C_1n^{-c_1}$, then
%
%
\begin{equation}
\label{eqmainresultcorexactstrong} \sup_{f\in\mathcal{F}}\bigl\llvert \Pr_f(f\in
\widetilde{\mathcal {C}}_{n})-(1-\alpha)\bigr\rrvert \leq
Cn^{-c}.
\end{equation}
Here $c$ and $C$ are constants depending only on $\alpha$, $c_1$ and $C_1$.
\end{corollary}

%

%

\subsection{Verifying Conditions \texorpdfstring{\protect\ref{con1}--\protect\ref{con4}}{H1-H4} for
confidence bands constructed using common density estimators via Gaussian multiplier bootstrap}
We now argue that when $\hat{c}_n(\alpha)$ is constructed via
Gaussian multiplier bootstrap, Conditions~\ref{con1}--\ref{con4}
hold for common density estimators---specifically, both for convolution
and for projection kernel density estimators under mild assumptions on
the kernel function.

Let $\{K_l\}_{l\in\mathcal{L}_n}$ be a family of kernel functions
where $K_l\dvtx  \mathbb{R}^d\times\mathbb{R}^d \to\mathbb{R}$ and $l$
is a smoothing parameter. We consider kernel density estimators of the form
%
%
\begin{equation}
\label{eqkerneldensityestimator} \hat{f}_{n}(x,l):= \mathbb{E}_n\bigl[
K_{l} (X_i,x) \bigr] = \frac{1}{n} \sum
_{i=1}^{n} K_{l} (X_{i},x),
\end{equation}
where $x\in\mathcal{X}$ and $l\in\mathcal{L}_n$.
The variance of $\sqrt{n}\hat{f}_n(x,l)$ is given by
\[
\sigma^{2}_{n,f}(x,l): = \Ep_f
\bigl[K_{l} (X_{1},x)^{2}\bigr] - \bigl(
\Ep_f\bigl[K_{l} (X_{1},x)\bigr]
\bigr)^{2}.
\]
We estimate $\sigma^{2}_{n,f}(x,l)$ by
%
%
\begin{equation}
\label{eqsigmaestimator} \hat{\sigma}_{n}^{2}(x,l):=
\frac{1}{n} \sum_{i=1}^{n}
K_{l}(X_{i},x)^{2} - \hat{f}_{n}(x,l)^{2}.
\end{equation}
This is a sample analogue estimator.

\subsubsection*{Examples}
Our general theory covers a wide class of kernel functions, such as
convolution, wavelet projection and nonwavelet projection kernels.
\begin{longlist}[(iii)]
\item[(i)] \textit{Convolution kernel}. Consider a function $K\dvtx \mathbb
{R}\rightarrow\mathbb{R}$. Let $\mathcal{L}_n\subset(0,\infty)$.
Then\vadjust{\goodbreak} for $x=(x_1,\ldots,x_d)'\in\mathbb{R}^d$, $y=(y_1,\ldots,y_d)'\in\mathbb{R}^d$ and $l\in\mathcal{L}_n$, the convolution
kernel function is defined by
%
%
\begin{equation}
\label{eqconvolutionkernel} K_l(y,x):=2^{ld}\prod
_{1\leq m\leq d}K \bigl(2^l(y_m-x_m)
\bigr).
\end{equation}
Here $2^{-l}$ is the bandwidth parameter.

\item[(ii)] \textit{Wavelet projection kernel}. 
Consider a father wavelet $\phi$, that is, a function $\phi$ such
that: (a) $\{\phi(\cdot-k)\dvtx k\in\mathbb{Z}\}$ is an orthonormal
system in $L_2(\mathbb{R})$, (b) the spaces $V_j=\{\sum_{k}c_k\phi
(2^jx-k)\dvtx \sum_{k}c_k^2<\infty\}$, $j=0,1,2,\ldots,$ are nested in the
sense that $V_j\subset V_{j^{\prime}}$ whenever $j\leq j^{\prime}$
and (c) $\bigcup_{j\geq0}V_j$ is dense in $L_2(\mathbb{R})$. Let
$\mathcal{L}_n\subset\mathbb{N}$. Then for $x=(x_1,\ldots,x_d)'\in
\mathbb{R}^d$, $y=(y_1,\ldots,y_d)'\in\mathbb{R}^d$, and $l\in
\mathcal{L}_n$, the wavelet projection kernel function is defined by
%
%
\begin{equation}
\label{eqwaveletkernel} K_l(y,x):=2^{ld}\sum
_{k_1,\ldots,k_d\in\mathbb{Z}}\prod_{1\leq
m\leq d}\phi
\bigl(2^ly_m-k_m\bigr)\prod
_{1\leq m\leq d}\phi\bigl(2^lx_m-k_m
\bigr).
\end{equation}
Here $l$ is the resolution level. We refer to \cite{Daubechies92} and
\cite{Hardle} as basic references on wavelet theory.

\item[(iii)] \textit{Nonwavelet projection kernel}. Let $\{\varphi
_j\dvtx j=1,\ldots, \infty\}$ be an orthonormal basis of $L_2(\mathcal
{X})$, the space of square integrable (with respect to Lebesgue
measure) functions on $\mathcal{X}$. Let $\mathcal{L}_n\subset
(0,\infty)$. Then for $x=(x_1,\ldots,x_d)'\in\mathbb{R}^d$,
$y=(y_1,\ldots,y_d)'\in\mathbb{R}^d$ and $l\in\mathcal{L}_n$, the
nonwavelet projection kernel function is defined by
%
%
\begin{equation}
\label{eqnon-waveletkernel} K_l(y,x):=\sum_{j=1}^{\lfloor2^{l d}\rfloor}
\varphi_j(y)\varphi_j(x),
\end{equation}
where $\lfloor a\rfloor$ is the largest integer that is smaller than
or equal to $a$.
Here $\lfloor2^{l d}\rfloor$ is the number of series (basis) terms
used in the estimation. When $d=1$ and $\mathcal{X}=[-1,1]$, examples
of orthonormal bases are Fourier basis
%
%
\begin{equation}
\label{eqfourierseries} \bigl\{1,\cos(\pi x),\cos(2\pi x),\ldots\bigr\}
\end{equation}
and Legendre polynomial basis
%
%
\begin{equation}
\label{eqpolynomialbasis} \bigl\{1,(3/2)^{1/2}x,(5/8)^{1/2}
\bigl(3x^2-1\bigr),\ldots\bigr\}.
\end{equation}
When $d>1$ and $\mathcal{X}=[-1,1]^d$, one can take tensor products of
bases for $d=1$.
\end{longlist}

We assume that the critical value $\hat{c}_n(\alpha)$ is obtained via
the multiplier bootstrap method:

%
\begin{algorithm}[(Gaussian multiplier bootstrap)]\label{al1}
Let $\xi_1,\ldots,\xi_n$ be independent $N(0,1)$ random variables
that are independent of the data $X_1^n:=\{X_1,\ldots,X_n\}$. Let $\xi
_1^n:=\{\xi_1,\ldots,\xi_n\}$. For all $x\in\mathcal{X}$ and $l\in
\mathcal{L}_n$, define\vadjust{\goodbreak} a Gaussian multiplier process
%
%
\begin{equation}
\label{eqbootstrapprocess} \hat{\mathbb{G}}_n(x,l):=\hat{\mathbb{G}}_n
\bigl(X_1^n,\xi _1^n\bigr)
(x,l):=\frac{1}{\sqrt{n}}\sum_{i=1}^n
\xi_i\frac
{K_l(X_i,x)-\hat{f}_n(x,l)}{\hat{\sigma}_n(x,l)}.
\end{equation}
Then the estimated critical value $\hat{c}_n(\alpha)$ is defined as
\[
\hat{c}_n(\alpha) = \mbox{conditional $(1-\alpha)$-quantile of $
\Vert\hat{\mathbb{G}}_n\Vert_{\mathcal{V}_n}$ given
$X_1^n$}.
\]
\end{algorithm}
Gaussian multiplier bootstrap is a special case of a more general
exchangeable bootstrap; see, for example, \cite{Wellner93}. We refer
the reader to \cite{GZ84} for the first systematic use of the Gaussian
multipliers and to \cite{LT88} and \cite{GZ90} for conditional
multiplier central limit theorems in the Donsker setting.

Let
\[
\mathcal{K}_{n,f}:= \biggl\{\frac{K_l(\cdot,x)}{\sigma
_{n,f}(x,l)}\dvtx (x,l)\in\mathcal{X}\times
\mathcal{L}_n \biggr\}
\]
denote the class of studentized kernel functions,
and define
\[
\sigma_n=\sup_{f\in\mathcal{F}}\sup_{g\in\mathcal{K}_{n,f}}
\bigl(\Ep_f\bigl[g(X_1)^2\bigr]
\bigr)^{1/2}.
\]
Note that $\sigma_n\geq1$.

For a given class $\mathcal{G}$ of measurable functions on a
probability space $(S,\mathcal{S},Q)$ and $\epsilon> 0$, the
$\epsilon$-covering number of $\mathcal{G}$ with respect to the
$L_{2}(Q)$-semimetric is denoted by $N(\mathcal{G}, L_{2}(Q), \epsilon
)$; see Chapter~2 of \cite{VW96} on details of covering numbers.
We will use the following definition of VC type classes:
%
\begin{definition}[(VC type class)]
Let $\mathcal{G}$ be a class of measurable functions on a measurable
space $(S,\mathcal{S})$, and let $b>0$, $a\geq e$ and $v\geq1$ be
some constants. Then the class $\mathcal{G}$ is called VC$(b,a,v)$
type class if it is uniformly bounded in absolute value by $b$ (i.e.,
$\sup_{g \in\mathcal{G}} \| g \|_{\infty} \leq b$), and the
covering numbers of $\mathcal{G}$ satisfy
\[
\sup_Q N\bigl(\mathcal{G},L_2(Q),b\tau\bigr)
\leq(a/\tau)^v,\qquad 0<\tau<1,
\]
where the supremum is taken over all finitely discrete probability
measures $Q$ on~$(S,\mathcal{S})$.
\end{definition}

Then we will assume the following condition.

\renewcommand{\theconditionvc}{VC}
%
\begin{conditionvc}\label{conVCtype}
There exist sequences $b_n>0$, $a_n\geq e$ and $v_n\geq1$ such that
for every $f \in\mathcal{F}$, the class $\mathcal{K}_{n,f}$ is
VC$(b_n,a_n,v_n)$ type and pointwise measurable.
\end{conditionvc}
We refer to Chapter~2.3 of \cite{VW96} for the definition of pointwise
measurable classes of functions.
We note that Condition~\ref{conVCtype} is a mild assumption, which we verify for
common constructions in Appendix~F (as a part of
proving results for the next section; see Comment~3.5 below); see also
Appendix~I (supplemental material \cite{CCK2014}).

For some sufficiently large absolute constant $A$, take
\[
K_n:=Av_n \bigl(\log n\vee\log(a_nb_n/
\sigma_n) \bigr).
\]
We will assume without loss of generality that $K_n\geq1$ for all $n$.
The following theorem verifies Conditions~\ref{con1}--\ref{con4}
with so defined $\hat{\sigma}_{n}^{2}(x,l)$ and $\hat{c}_{n}(\alpha
)$ under Condition~\ref{conVCtype}, using the critical values constructed via
Algorithm~\ref{al1}.

%
\begin{theorem}[(Conditions~\ref{con1}--\ref{con4} hold for our construction)]\label{propverifyingconditions1-3}
Suppose that Condition~\ref{conVCtype} is satisfied and there exist constants
$c_{2},C_{2} > 0$ such that $b_n^2\sigma_n^4K_n^4/n\leq C_2n^{-c_2}$.
Then Conditions~\ref{con1}--\ref{con4}, including both Conditions~\ref{con2}\textup{(a)} and~\ref{con2}\textup{(b)}, hold with some constants $c_1, C_1>0$
that depend only on $c_2, C_2$.
\end{theorem}

\section{Honest and adaptive confidence bands in H\"{o}lder classes}\label{secverifyingconditions}
In this section, we study the problem of constructing honest adaptive
confidence bands in H\"{o}lder smoothness classes.
Recall that for $t, L > 0$, the\vadjust{\goodbreak} H\"{o}lder ball of densities with
radius $L$ and smoothness level $t$ is defined by
\begin{eqnarray*}
\Sigma(t,L)&:=& \biggl\{ f\dvtx  \mathbb{R}^{d} \to\mathbb{R}\dvtx  f\mbox{ is
a $\lfloor t \rfloor$-times continuously differentiable density},
\\
&&\hspace*{5pt} \bigl\| D^{\alpha} f \bigr\|_{\infty} \leq L,\ \forall| \alpha| \leq\lfloor t
\rfloor,
\sup_{x \neq y} \frac{| D^{\alpha} f (x) - D^{\alpha
}f(y)|}{|x-y|^{t-\lfloor t \rfloor}} \leq L,
\\[-6pt]
&&\hspace*{247pt} \forall|\alpha| = \lfloor t \rfloor \biggr\},
\end{eqnarray*}
where $\lfloor t\rfloor$ denotes the largest integer smaller than $t$,
and for a multi-index $\alpha= (\alpha_{1},\ldots,\alpha_{d})$ with
$| \alpha| = \alpha_{1}+\cdots+\alpha_{d}$, $D^{\alpha}f (x):=
\partial^{|\alpha|} f(x)/\partial x_{1}^{\alpha_{1}} \cdots\partial
x_{d}^{\alpha_{d}}$; see, for example, \cite{Tsybakov09}. We assume
that for some $0<\underline{t}\leq\bar{t}<\infty$ and $L \geq1$,
%
%
\begin{equation}
\label{eqHolderclasses} \mathcal{F}\subset\bigcup_{t\in[\underline{t},\bar{t}]}
\Sigma(t,L),
\end{equation}
and consider the confidence band $\mathcal{C}_{n} = \{ \mathcal
{C}_{n}(x)\dvtx  x \in\mathcal{X} \}$ of the form (\ref{cband}) and
(\ref{eqsn}), where $\mathcal{X}$ is a (suitable) compact set in
$\mathbb{R}^d$.

We begin by stating our assumptions. First, we restrict attention to
kernel density estimators $\hat{f}_n$ based on either convolution or
wavelet projection kernel functions. Let $r$ be an integer such that $r
\geq2$ and $r > \bar{t}$.

\renewcommand{\theconditionlow}{L\arabic{conditionlow}}
%
\begin{conditionlow}[(Density estimator)]\label{condensityestimator}
The density estimator $\hat{f}_n$ is either a convolution or wavelet
projection kernel density estimator defined in (\ref
{eqkerneldensityestimator}), (\ref{eqconvolutionkernel})~and~(\ref
{eqwaveletkernel}). For convolution kernels, the function $K\dvtx  \mathbb
{R}\to\mathbb{R}$ has compact support and is of bounded variation,
and moreover is such that $\int K(s) \,ds = 1$ and $\int s^{j} K(s) \,dx =
0$ for $j=1,\ldots,r-1$. For wavelet projection kernels, the function
$\phi\dvtx  \mathbb{R}\to\mathbb{R}$ is either a compactly supported
father wavelet of regularity $r-1$ [i.e., $\phi$ is $(r-1)$-times
continuously differentiable], or a Battle--Lemari\'{e} wavelet of
regularity $r-1$.
\end{conditionlow}
The assumptions stated in Condition~\ref{condensityestimator} are
commonly used in the literature. See \cite{GG02} for a more general
class of convolution kernel functions that would suffice for our results.
Details on compactly supported and Battle--Lemari\'{e} wavelets can be
found in Chapters~6 and 5.4 of \cite{Daubechies92}, respectively.

It is known that if the function class $\mathcal{F}$ is sufficiently
large [e.g., if $\mathcal{F}=\Sigma(t,L)\cup\Sigma(t',L)$ for
$t'>t$], the construction of honest adaptive confidence bands is not
possible; see \cite{Low97}. Therefore, following \cite{GN10}, we will
restrict the function class $\mathcal{F}\subset\bigcup_{t\in
[\underline{t},\bar{t}]}\Sigma(t,L)$ in a suitable way, as follows:

%
\begin{conditionlow}[(Bias bounds)]\label{congineandnickl}
There exist constants $l_0, c_3, C_3 > 0$ such that for every $f\in
\mathcal{F}\subset\bigcup_{t\in[\underline{t},\bar{t}]}\Sigma
(t,L)$, there exists $t\in[\underline{t},\bar{t}]$ with
%
%
\begin{equation}
\label{eqbiasbounds} c_32^{-lt}\leq\sup_{x\in\mathcal{X}}\bigl|
\Ep_f\bigl[\hat {f}_n(x,l)\bigr]-f(x)\bigr|\leq
C_32^{-lt},
\end{equation}
for all $l\geq l_0$.\vadjust{\goodbreak}
\end{conditionlow}

This condition is inspired by the path-breaking work of \cite{GN10};
see also \cite{PT00}. It can be interpreted as the requirement that
the functions $f$ in the class $\mathcal{F}$ are ``self-similar'' in
the sense that their regularity remains the same at large and small
scales; see also \cite{B11a}. To put it differently,
``self-similarity'' could be understood as the requirement that the
bias of the kernel\vspace*{1pt} approximation to $f$ with bandwidth $2^{-l}$ remains
approximately proportional to $(2^{-l})^t$---that is, not much smaller
or not much bigger---for all small values of the bandwidth~$2^{-l}$.

It is useful to note that the upper bound in (\ref{eqbiasbounds})
holds for all $f\in\Sigma(t,L)$ (for sufficiently large $C_3$) under
Condition~\ref{condensityestimator}; see, for example, Theorem~9.3 in
\cite{Hardle}. In addition, Gin\'{e} and Nickl \cite{GN10} showed
that under Condition~\ref{condensityestimator}, the restriction due
to the lower bound in (\ref{eqbiasbounds}) is weak in the sense that
the set of elements of $\Sigma(t,L)$ for which the lower bound in
(\ref{eqbiasbounds}) does not hold is ``topologically small.''
Moreover, they showed that the minimax optimal rate of convergence in
the sup-norm over $\Sigma(t,L)$ coincides with that over the set of
elements of $\Sigma(t,L)$ for which Condition~\ref{congineandnickl}
holds. We refer to \cite{GN10} for a detailed and deep discussion of
these conditions and results.

We also note that, depending on the problem, construction of honest
adaptive confidence bands is often possible under somewhat weaker
conditions than that in Condition~\ref{congineandnickl}. For example, if we are
interested in the function class $\Sigma(t,L)\cup\Sigma(t',L)$ for
some $t'>t$, Hoffman and Nickl \cite{HN11} showed that it is necessary
and sufficient to exclude functions $\Sigma(t,L)\setminus\Sigma
(t,L,\rho_n)$ where $\Sigma(t,L,\rho_n)=\{f\in\Sigma(t,L)\dvtx \inf_{g\in\Sigma(t',L)}\Vert g-f\Vert_{\infty}\geq\rho_n\}$ and where
$\rho_n>0$ is allowed to converge to zero as $n$ increases but
sufficiently slowly. If we are interested in the function class $\bigcup_{t\in[\underline{t},\bar{t}]}\Sigma(t,L)$, Bull \cite{B11a}
showed that (essentially) necessary and sufficient condition can be
written in the form of the bound from below on the rate with which
wavelet coefficients of the density $f$ are allowed to decrease. Here
we prefer to work with Condition~\ref{congineandnickl} directly
because it is directly related to the properties of the estimator $\hat
{f}_n$ and does not require any further specifications of the function
class $\mathcal{F}$.

In order to introduce the next condition, we need to observe that under
Condition~\ref{congineandnickl}, for every $f\in\mathcal{F}$, there
exists a \textit{unique} $t\in[\underline{t},\bar{t}]$ satisfying
(\ref{eqbiasbounds}); indeed, if $t_1<t_2$, then for any $c,C>0$,
there exists $\bar{l}$ such that $C2^{-lt_2}<c2^{-lt_1}$ for all
$l\geq\bar{l}$, so that for each $f\in\mathcal{F}$ condition (\ref
{eqbiasbounds}) can hold for all $l\geq l_0$ for at most one value of
$t$. This defines the map
%
%
\begin{equation}
\label{eqdefinet} t\dvtx \mathcal{F} \to[\underline{t},\bar{t}],\qquad f \mapsto t(f).
\end{equation}
%
The next condition states our assumptions on the candidate set
$\mathcal{L}_n$ of the values of the smoothing parameter:
%
\begin{conditionlow}[(Candidate set)]\label{concandidateset}
There exist constants $c_{4},C_{4} >0$ such that for every $f\in
\mathcal{F}$, there exists $l\in\mathcal{L}_n$ with
%
%
\begin{equation}
\label{eql0star} \biggl(\frac{c_4\log n}{n} \biggr)^{1/(2t(f)+d)}
\leq2^{-l}\leq \biggl(\frac{C_4\log n}{n} \biggr)^{1/(2t(f)+d)},
\end{equation}
for the map $t\dvtx  f \mapsto t(f)$ defined in (\ref{eqdefinet}). In
addition, the candidate set is $\mathcal{L}_n=[l_{\min,n},l_{\max,n}]\cap\mathbb{N}$.
\end{conditionlow}
This condition thus ensures via (\ref{eql0star}) that the candidate
set $\mathcal{L}_n$ contains an appropriate value of the smoothing
parameter that leads to the optimal rate of convergence for every
density $f\in\mathcal{F}$.

Finally, we will make the following mild condition:
%
\begin{conditionlow}[(Density bounds)]\label{condensitybounds}
There exist constants $\delta, \underline{f},\bar{f} > 0$ such that
for all $f\in\mathcal{F}$,
%
%
\begin{equation}
\label{eqdensityconstraint} f(x)\geq\underline{f}\qquad\mbox{for all }x\in
\mathcal{X}^\delta \quad\mbox{and}\quad f(x)\leq\bar{f}\qquad\mbox{for all
}x\in \mathbb{R}^d,
\end{equation}
where $\mathcal{X}^\delta$ is the $\delta$-enlargement of $\mathcal
{X}$, that is, \mbox{$\mathcal{X}^{\delta} = \{ x \in\mathbb{R}^{d}\dvtx
\inf_{y \in\mathcal{X}} | x - y | \leq\delta\}$}.
\end{conditionlow}

We now discuss how we choose various parameters in the confidence band
$\mathcal{C}_{n}$.
In the previous section, we have shown how to obtain honest confidence
bands as long as we can control the bias $\Delta_{n,f}(\hat{l}_n)$
appropriately. So to construct honest adaptive confidence bands, we
seek a method to choose the smoothing parameter $\hat{l}_n\in\mathcal
{L}_n$ so that the bias $\Delta_{n,f}(\hat{l}_n)$ can be controlled,
and at the same time, the confidence band $\mathcal{C}_{n}$ is adaptive.

Let $\widebar{\mathcal{V}}_n:=\{(x,l,l')\dvtx x\in\mathcal{X},l,l'\in
\mathcal{L}_n,l<l'\}$, and for $(x,l,l')\in\widebar{\mathcal{V}}_n$, denote
\[
\tilde{\sigma}_n\bigl(x,l,l'\bigr):= \Biggl(
\frac{1}{n}\sum_{i=1}^n
\bigl(K_l(X_i,x)-K_{l'}(X_i,x)
\bigr)^2-\bigl(\hat{f}_n(x,l)-\hat {f}_n
\bigl(x,l'\bigr)\bigr)^2 \Biggr)^{1/2}.
\]
Also, for some small $c_{\sigma}>0$, let
\[
\hat{\sigma}_n\bigl(x,l,l'\bigr):=
\bigl(c_\sigma\hat{\sigma}_n\bigl(x,l'\bigr)
\bigr)\vee \tilde{\sigma}_n\bigl(x,l,l'\bigr)
\]
denote the truncated version of $\tilde{\sigma}_n(x,l,l')$. In
practice, we suggest setting $c_\sigma=0.5(1-2^{-d/2})$ (the constant
$c_\sigma$ is chosen so that with probability approaching one, $\hat
{\sigma}_n(x,l,l')=\tilde{\sigma}_n(x,l,l')$ for all
$(x,l,l')\in\widebar{\mathcal{V}}_n$ for convolution kernel estimators,
and for all $(x,l,l')\in\widebar{\mathcal{V}}_n$ with $l\leq l'-s$ for
wavelet projection kernel estimators where~$s$ is some constant; see
Lemmas~F.2 and~F.4 in the supplemental material \cite{CCK2014}).

There exist several techniques in the literature to construct $\hat
{l}_{n}$ so that $\Delta_{n,f}(\hat{l}_n)$ can be controlled and the
confidence band $\mathcal{C}_n$ is adaptive; see, for example, \cite
{M07} for a thorough introduction. One of the most important such
techniques is the Lepski method; see \cite{Lepski91} for a detailed
explanation of the method. In this paper, we introduce a new
implementation of the Lepski method, which we refer to as a multiplier
bootstrap implementation of the Lepski method.

%
\begin{algorithm}[(Multiplier bootstrap implementation of the Lepski method)]
Let $\gamma_n$ be a sequence of positive numbers converging to zero.
Let $\xi_1,\ldots,\xi_n$ be independent $N(0,1)$ random variables
that are independent of the data $X_1^n:=\{X_1,\ldots,X_n\}$. Let $\xi
_1^n:=(\xi_1,\ldots,\xi_n)$. For all $(x,l,l')\in\widebar{\mathcal
{V}}_n$, define a Gaussian multiplier process
\begin{eqnarray*}
\widetilde{\mathbb{G}}_n\bigl(x,l,l'\bigr)&:=&
\widetilde{\mathbb {G}}_n\bigl(X_1^n,
\xi_1^n\bigr) \bigl(x,l,l'\bigr)
\\
&:=&\frac{1}{\sqrt{n}}\sum_{i=1}^n
\xi_i\frac
{(K_l(X_i,x)-K_{l'}(X_i,x))-(\hat{f}_n(x,l)-\hat{f}_n(x,l'))}{\hat
{\sigma}_n(x,l,l')}.
\end{eqnarray*}
Also, define
\[
\tilde{c}_n(\gamma_n)=\mbox{conditional $(1-
\gamma_n)$-quantile of $\|\widetilde{\mathbb{G}}_n
\|_{\widebar{\mathcal{V}}_n}$ given $X_1^n$}. %
\]
%
Morever, for all $l\in\mathcal{L}_n$, let
\[
\mathcal{L}_{n,l}:=\bigl\{l'\in\mathcal{L}_n\dvtx l'>l
\bigr\}.
\]
Finally, for some constant $q>1$, which is independent of $n$, define a
Lepski-type estimator
%
%
\begin{equation}
\label{eqlepski-typeestimator} \hat{l}_n:=\inf \biggl\{l\in\mathcal{L}_n\dvtx
\sup_{l'\in\mathcal
{L}_{n,l}}\sup_{x\in\mathcal{X}}\frac{\sqrt{n}|\hat
{f}_n(x,l)-\hat{f}_n(x,l')|}{\hat{\sigma}_n(x,l,l')}\leq q
\tilde{c}_n(\gamma_n) \biggr\}.
\end{equation}
\end{algorithm}

%
Once we have $\hat{l}_n$, to define the confidence band $\mathcal
{C}_{n}$, we need to specify $\hat{\sigma}_{n}(x,l)$, $\hat
{c}_n(\alpha)$ and $c_n'$. We assume that $\hat{\sigma}_n(x,l)$ is
obtained via (\ref{eqsigmaestimator}) and $\hat{c}_n(\alpha)$ via
Algorithm~\ref{al1}. To specify $c_n'$, let $u_n'$ be a sequence of positive
numbers such that $u_n'$ is sufficiently large for large $n$.
Specifically, for large $n$, $u_n'$ is assumed to be larger than some
constant $C(\mathcal{F})$ depending only on the function class
$\mathcal{F}$. Set
\[
c_n':=u_n'
\tilde{c}_n(\gamma_n).
\]

The following theorem shows that the confidence band $\mathcal{C}_{n}$
defined in this way is honest and adaptive for $\mathcal{F}$:

%
\begin{theorem}[(Honest and adaptive confidence bands via our method)]\label{propverifyingconditionH4}
Suppose that Conditions~\ref{condensityestimator}--\ref{condensitybounds} are satisfied. In addition, suppose that there exist
constants $c_{5},C_{5} > 0$ such that: \textup{(i)} $2^{l_{\max,n}d}(\log^4
n)/n\leq C_5n^{-c_5}$,\break  \textup{(ii)} $l_{\min,n}\geq c_5\log n$, \textup{(iii)}~$\gamma
_n\leq C_5n^{-c_5}$, \textup{(iv)} $|\log\gamma_n|\leq C_5\log n$, \textup{(v)}
$u_n'\geq C(\mathcal{F})$ and \textup{(vi)} $u_n'\leq C_5\log n$. Then
Conditions \ref{con1}--\ref{con6} in Section~\ref{secanalysisunderhighlevelconditions}
and (\ref{eqresultstrong}) in Theorem~\ref{prop1} hold and
%
%
\begin{equation}
\label{eqintervallength} \sup_{f \in\mathcal{F}} \Pr_{f} \Bigl( \sup
_{x \in\mathcal{X}} \lambda\bigl(\mathcal{C}_{n}(x)\bigr) > C
\bigl(1+u_n'\bigr) r_n\bigl(t(f)\bigr)
\Bigr)\leq Cn^{-c},
\end{equation}
where $\lambda(\cdot)$ denotes the Lebesgue measure on $\mathbb{R}$
and $r_n(t):=(\log n/n)^{t/(2t+d)}$. Here the constants $c, C > 0$
depend only on $c_{5},C_{5}$, the constants that appear in Conditions
\ref{condensityestimator}--\ref{condensitybounds}, $c_\sigma$,
$\alpha$ and the function $K$ (when convolution kernels are used) or
the father wavelet $\phi$ (when wavelet projection kernels are used). Moreover,
%
%
\begin{equation}
\label{eqadaptivity} \sup_{f\in\mathcal{F}\cap\Sigma(t,L)}\Pr_{f} \Bigl( \sup
_{x \in
\mathcal{X}} \lambda\bigl(\mathcal{C}_{n}(x)\bigr) > C
\bigl(1+u_n'\bigr) r_n(t) \Bigr)\leq
Cn^{-c},
\end{equation}
with the same constants $c,C$ as those in (\ref{eqintervallength}).
\end{theorem}

\begin{appendix}\label{app}
\section{Coupling inequalities for suprema of empirical and related processes}\label{seccouplinginequalities}
The purpose of this section is to provide two coupling inequalities
based on Slepian--Stein methods that are useful for the analysis of
uniform confidence bands. The first inequality is concerned with
suprema of empirical processes and is proven in Corollary~2.2 in \cite
{CCK2012}. The second inequality is new, is concerned with suprema of
Gaussian multiplier processes, and will be obtained from a Gaussian
comparison theorem derived in \cite{CCK2012b}.



Let $X_1,\ldots,X_n$ be i.i.d. random variables taking values in a
measurable space $(S,\mathcal{S})$. Let $\mathcal{G}$ be a
pointwise-measurable VC$(b,a,v)$ type function class for some $b>0$,
$a\geq e$, and $v\geq1$ (the definition of VC type classes is given in
Section~\ref{secanalysisunderhighlevelconditions}). Let $\sigma^2 >
0$ be any constant such that $\sup_{g\in\mathcal{G}}\Ep
[g(X_1)^2]\leq\sigma^2\leq b^2$.
Define the empirical process
\[
\mathbb{G}_n(g):=\frac{1}{\sqrt{n}}\sum_{i=1}^n
\bigl(g(X_i)-\Ep \bigl[g(X_1)\bigr] \bigr),\qquad g\in
\mathcal{G},
\]
and let
\[
W_n:=\Vert\mathbb{G}_n\Vert_{\mathcal{G}}:=\sup
_{g\in\mathcal
{G}}\bigl|\mathbb{G}_n(g)\bigr|
\]
denote the supremum of the empirical process. Note that $W_n$ is a
well-defined random variable since $\mathcal{G}$ is assumed to be
pointwise-measurable.
Let $B=\{B(g)\dvtx g\in\mathcal{G}\}$ be a tight Gaussian random element
in $\ell^{\infty}(\mathcal{F})$ with mean zero and covariance function
\[
\Ep\bigl[B(g_1)B(g_2)\bigr]=\Ep\bigl[g_1(X_1)g_2(X_1)
\bigr]-\Ep\bigl[g_1(X_1)\bigr]\Ep\bigl[g_2(X_1)
\bigr],
\]
for all $g_1,g_2\in\mathcal{G}$. It is well known that such a process
exists under the VC type assumption; see \cite{VW96}, pages 100--101.
Finally, for some sufficiently large absolute constant $A$, let
\[
K_n:=Av\bigl(\log n\vee\log(ab/\sigma)\bigr).
\]
In particular, we will assume that $K_n\geq1$.
The following theorem shows that $W_n$ can be well approximated by the
supremum of the corresponding Gaussian process~$B$ under mild
conditions on $b$, $\sigma$ and $K_n$. The proof of this theorem can
be found in Corollary~2.2 in \cite{CCK2012}.
%
\begin{theorem}[(Slepian--Stein type coupling for suprema of empirical processes)]\label{thmcouplinginequality1}
Consider the setting specified above.
Then for every $\gamma\in(0,1)$ one can construct on an enriched
probability space a random variable $W^0$ such that: \textup{(i)}~$W^0\stackrel
{d}{=}\Vert B\Vert_{\mathcal{G}}$ and \textup{(ii)}
\begin{eqnarray*}
&& \Pr \biggl(\bigl|W_n-W^0\bigr|>\frac{bK_n}{(\gamma n)^{1/2}}+
\frac{(b\sigma
)^{1/2}K_n^{3/4}}{\gamma^{1/2}n^{1/4}}+\frac{b^{1/3}\sigma
^{2/3}K_n^{2/3}}{\gamma^{1/3}n^{1/6}} \biggr)
\\
&&\qquad\leq A' \biggl(\gamma+\frac{\log n}{n} \biggr),
\end{eqnarray*}
where $A'$ is an absolute constant.
\end{theorem}
%
%

Let $\xi_1,\ldots,\xi_n$ be independent $N(0,1)$ random variables
independent of $X_1^n:=\{X_1,\ldots,X_n\}$, and let $\xi_1^n:=\{\xi
_1,\ldots,\xi_n\}$. We assume that random variables $X_1,\ldots,X_n,\xi_1,\ldots,\xi_n$ are\vadjust{\goodbreak} defined as coordinate projections from
the product probability space.
Define the Gaussian multiplier process
\[
\widetilde{\mathbb{G}}_n(g):=\widetilde{\mathbb{G}}_n
\bigl(X_1^n,\xi _1^n\bigr) (g):=
\frac{1}{\sqrt{n}}\sum_{i=1}^n
\xi_i\bigl(g(X_i)-\mathbb {E}_n
\bigl[g(X_i)\bigr]\bigr),\qquad g\in\mathcal{G},
\]
and\vspace*{1pt} for $x_1^n\in S^n$, let $\widetilde{W}_n(x_1^n):=\Vert\widetilde{\mathbb
{G}}_n(x_1^n,\xi_1^n)\Vert_{\mathcal{G}}$ denote the supremum of
this process calculated for fixed $X_1^n=x_1^n$. Note that $\widetilde
{W}_n(x_1^n)$ is a well-defined random variable.
In addition, let
\[
\psi_n:=\sqrt{\frac{\sigma^2 K_n}{n}}+ \biggl(\frac{b^2\sigma^2
K_n^3}{n}
\biggr)^{1/4}\quad\mbox{and}\quad \gamma_n(\delta):=
\frac{1}{\delta} \biggl(\frac{b^2\sigma^2
K_n^3}{n} \biggr)^{1/4}+
\frac{1}{n}.
\]
The following theorem shows that $\widetilde{W}_n(X_1^n)$ can be well
approximated with high probability by the supremum of the Gaussian
process $B$ under mild conditions on~$b$, $\sigma$ and $K_n$. The
proof of this theorem can be found in the supplemental material~\cite{CCK2014}.
%
\begin{theorem}[(Slepian--Stein type coupling for suprema of conditional multiplier processes)]\label{thmcouplinginequality2}
Consider the setting specified above. Suppose that $b^2K_n\leq n\sigma^2$.
Then for every $\delta>0$, there exists a set $S_{n,0}\in\mathcal
{S}^n$ such that $\Pr( X_{1}^{n} \in S_{n,0})\geq1-3/n$ and for every
$x_1^n\in S_{n,0}$ one can construct on an enriched probability space a
random variable $W^0$ such that: \textup{(i)} $W^0\stackrel{d}{=}\Vert B\Vert
_{\mathcal{G}}$ and \textup{(ii)}
\[
\Pr\bigl(\bigl|\widetilde{W}_n\bigl(x_1^n
\bigr)-W^0\bigr|>(\psi_n+\delta)\bigr)\leq A''
\gamma _n(\delta),
\]
where $A''$ is an absolute constant.
\end{theorem}
%
%

\section{Some technical tools}\label{sectechnicaltools}

%
\begin{theorem}\label{thmtalagrandinequality}
Let $\xi_{1},\ldots,\xi_{n}$ be i.i.d. random variables taking values
in a measurable space $(S,\mathcal{S})$.
Suppose that $\mathcal{G}$ is a nonempty, pointwise measurable class
of functions on $S$ uniformly bounded by a constant $b$ such that there
exist constants $a \geq e$ and $v > 1$ with
$\sup_{Q}N(\mathcal{G},L_{2}(Q),b\epsilon) \leq(a/\epsilon)^{v}$
for all $0 < \epsilon\leq1$. Let $\sigma^{2}$ be a constant such
that $\sup_{g \in \mathcal{G}} \var(g) \leq\sigma^{2} \leq
b^{2}$. If $b^{2} v\log(ab/\sigma) \leq n \sigma^{2}$,
then for all $t \leq n\sigma^2/b^2$,
\[
\Pr \Biggl[ \sup_{g \in\mathcal{G}} \Biggl\llvert \sum
_{i=1}^{n} \bigl\{ g(\xi _{i}) - \Ep
\bigl[g(\xi_{1})\bigr] \bigr\} \Biggr\rrvert > A \sqrt{ n
\sigma^{2} \biggl\{ t \vee \biggl(v \log\frac{ab}{\sigma} \biggr)
\biggr\} } \Biggr] \leq e^{-t},
\]
where $A > 0$ is an absolute constant.
\end{theorem}
\begin{pf}
This version of Talagrand's inequality follows from Theorem~3 in~\cite{M00} combined with a bound on expected values of suprema of empirical
processes derived in \cite{GG01}. See also \cite{T96} for the
original version of Talagrand's inequality.
\end{pf}

Proofs of the following two lemmas can be found in the supplemental material~\cite{CCK2014}.

%
\begin{lemma}\label{lemcriticalvaluebound}
Let $Y:=\{Y(t)\dvtx t\in T\}$ be a separable, centered Gaussian process such
that $\Ep[Y(t)^2]\leq1$ for all $t\in T$. Let $c(\alpha)$ denote the
$(1-\alpha)$-quantile of $\Vert Y\Vert_T$. Assume that $\Ep[\Vert
Y\Vert_T]<\infty$. Then $c(\alpha) \leq\Ep[\Vert Y\Vert_T] +\break
\sqrt{2 | \log \alpha|}$ and $c(\alpha) \leq M(\Vert Y\Vert_T) +
\sqrt{2 | \log \alpha|}$ for all $\alpha\in(0,1)$ where $M(\Vert
Y\Vert_T)$ is the median of $\Vert Y\Vert_T$.
\end{lemma}

%
\begin{lemma}
\label{lemvcproperties}
Let $\mathcal{G}_1$ and $\mathcal{G}_2$ be VC$(b_1,a_1,v_1)$ and
VC$(b_2,a_2,v_2)$ type classes, respectively, on a measurable space
$(S,\mathcal{S})$. Let $a=(a_1^{v_1}a_2^{v_2})^{1/(v_1+v_2)}$. Then:
\textup{(i)}~$\mathcal{G}_1\cdot\mathcal{G}_2=\{g_1\cdot g_2\dvtx g_1\in\mathcal
{G}_1,g_2\in\mathcal{G}_2\}$ is VC$(b_1b_2,2a,v_1+v_2)$ type class,
\textup{(ii)} $\mathcal{G}_1-\mathcal{G}_2=\{g_1-g_2\dvtx g_1\in\mathcal
{G}_1,g_2\in\mathcal{G}_2\}$ is VC$(b_1+b_2,a,v_1+v_2)$ type class
and \textup{(iii)} $\mathcal{G}_1^2=\{g_1^2\dvtx g_1\in\mathcal{G}_1\}$ is
VC$(b_1^2,2a_1,v_1)$ type class.
\end{lemma}

\section{Proofs for Section~\texorpdfstring{\lowercase{\protect\ref{secanti-concentration}}}{2}}\vspace*{-12pt}

\begin{pf*}{Proof of Theorem~\ref{thm2}}
The fact that $a(X) < \infty$ follows from Landau--Shepp--Fernique
theorem; see, for example, Lemma~2.2.5 in \cite{D99}. In addition,
since $\sup_{t \in T} X_{t}\geq X_{t_{0}}$ for any fixed $t_{0} \in
T$, $a(X) \geq\Ep[ X_{t_{0}}] = 0$. We now prove~(\ref
{eqanti-concentrationmain}).

Since the Gaussian process $X=(X_{t})_{t \in T}$ is separable, there
exists a sequence of finite subsets $T_{n} \subset T$ such that $Z_{n}:= \max_{t \in T_{n}} X_{t} \to\sup_{t \in T} X_{t} =: Z$ a.s. as $n
\to\infty$.
Fix any $x\in\mathbb{R}$. Since $|Z_{n} - x| \to|Z-x|$ a.s. and a.s.
convergence implies weak convergence, there exists an at most countable
subset $\mathcal{N}_{x}$ of $\mathbb{R}$ such that for all $\epsilon
\in\mathbb{R} \setminus\mathcal{N}_{x}$,
\[
\lim_{n \to\infty} \Pr\bigl(|Z_{n} - x| \leq\epsilon\bigr) = \Pr\bigl(| Z
- x | \leq\epsilon\bigr).
\]
But by Theorem~3 in \cite{CCK2012b},
\[
\Pr\bigl(|Z_n-x|\leq\epsilon\bigr)\leq4\epsilon\Bigl(\Ep\Bigl[\max
_{t\in
T_n}X_t\Bigr]+1\Bigr)\leq4\epsilon\bigl(a(X)+1
\bigr),
\]
for all $\epsilon\geq0$. Therefore,
%
%
\begin{equation}
\label{eqcountableset} \Pr\bigl(|Z-x|\leq\epsilon\bigr)\leq4\epsilon\bigl(a(X) + 1\bigr),
\end{equation}
for all $\epsilon\in\mathbb{R}\setminus\mathcal{N}_x$. By right
continuity of $\Pr(|Z-x|\leq\cdot)$, it follows that (\ref
{eqcountableset}) holds for all $\epsilon\geq0$. Since $x \in\mathbb
{R}$ is arbitrary, we obtain (\ref{eqanti-concentrationmain}).\vadjust{\goodbreak}
\end{pf*}

\begin{pf*}{Proof of Corollary \ref{cor1}}
In view of the proof of Theorem \ref{cor1}, it suffices to prove
the corollary in the case where $T$ is finite, but then the
corollary follows from Comment 5 in \cite{CCK2012b}.
\end{pf*}


\section{Proofs for Section~\texorpdfstring{\lowercase{\protect\ref{secanalysisunderhighlevelconditions}}}{3}}\vspace*{-12pt}
\begin{pf*}{Proof of Theorem~\ref{prop1}}
Pick any $f \in\mathcal{F}$.
By the triangle inequality, we have for any $x\in\mathcal{X}$,
\[
\frac{\sqrt{n}| \hat{f}_{n}(x,\hat{l}_{n}) - f(x) |}{\hat{\sigma
}_{n}(x,\hat{l}_{n})} \leq \bigl(\bigl|Z_{n,f}(x,\hat{l}_{n})\bigr|+
\Delta _{n,f}(\hat{l}_{n}) \bigr)\frac{\sigma_{n,f}(x,\hat{l}_n)}{\hat
{\sigma}_n(x,\hat{l}_n)},
\]
by which we have
\begin{eqnarray*}
\hspace*{-4pt}&& \Pr_{f}\bigl( f(x) \in\mathcal{C}_{n}(x),\ \forall x \in
\mathcal{X} \bigr)
\nonumber
\\
\hspace*{-4pt}&&\!\qquad \geq\Pr_{f} \bigl( \bigl| Z_{n,f}(x,\hat{l}_{n}) \bigr|+
\Delta_{n,f}(\hat {l}_n)\leq \bigl(\hat{c}_{n}(
\alpha)+c_n'\bigr)\hat{\sigma}_{n}(x,\hat
{l}_n)/\sigma_{n,f}(x,\hat{l}_n),\ \forall x \in
\mathcal{X} \bigr)
\nonumber
\\
\hspace*{-4pt}&&\!\qquad \geq\Pr_f\Bigl( \sup_{x
\in{\mathcal{X}}} \bigl|
Z_{n,f}(x,\hat{l}_{n}) \bigr|+\Delta_{n,f}(\hat
{l}_n)\leq \bigl(\hat{c}_{n}(\alpha)+c_n'
\bigr) (1-\epsilon_{3n})\Bigr)-\delta _{3n}
\\
\hspace*{-4pt}&&\!\qquad \geq\Pr_f\Bigl( \sup_{x\in\mathcal{X}}\bigl|
Z_{n,f}(x,\hat{l}_n)\bigr|\leq \hat{c}_{n}(\alpha) (1-
\epsilon_{3n})-c_n'\epsilon_{3n}
\Bigr)-\delta _{3n}-\delta_{4n}
\\
\hspace*{-4pt}&&\!\qquad \geq\Pr_f\bigl( \| Z_{n,f}\|_{\mathcal{V}_n}\leq
\hat{c}_{n}(\alpha ) (1-\epsilon_{3n})-c_n'
\epsilon_{3n}\bigr)-\delta_{3n}-\delta_{4n}
\\
\hspace*{-4pt}&&\!\qquad \geq\Pr_f\bigl( \| Z_{n,f}
\|_{\mathcal{V}_n}\leq \hat{c}_{n}(\alpha ) (1-\epsilon_{3n})-u_n
\epsilon_{3n}\sqrt{\log n}\bigr)-\delta _{3n}-
\delta_{4n}-\delta_{5n},
\end{eqnarray*}
where the third line follows from Condition~\ref{con4}, the fourth
line from Condition~\ref{con5}, the fifth line from the inequality
$\sup_{x\in\mathcal{X}}|Z_{n,f}(x,\hat{l}_n)|\leq\|Z_{n,f}\|
_{\mathcal{V}_n}$ and the sixth line from Condition~\ref{con6}.
Further, the probability in the last line above equals (recall that
$W_{n,f}=\|Z_{n,f}\|_{\mathcal{V}_n}$)
%
%
\begin{eqnarray}
&& \Pr_f\bigl( W_{n,f}\leq \hat{c}_{n}(\alpha)
(1-\epsilon _{3n})-u_n\epsilon_{3n}\sqrt{\log
n}\bigr)
\nonumber
\\
&&\qquad \geq\Pr_f\bigl(W_{n,f}\leq c_{n,f}(\alpha+
\tau_n) (1-\epsilon _{3n})-\epsilon_{2n}-u_n
\epsilon_{3n}\sqrt{\log n}\bigr)-\delta_{2n}, \label{eqprop1der5}
\end{eqnarray}
where (\ref{eqprop1der5}) follows from Condition~\ref{con3}. Now,
the probability in (\ref{eqprop1der5}) is bounded from below by
Condition~\ref{con1} by
%
%
\begin{eqnarray}
&&\Pr_f\bigl( W_{n,f}^0\leq c_{n,f}(
\alpha+\tau_n) (1-\epsilon _{3n})-\epsilon_{1n}-
\epsilon_{2n}-u_n\epsilon_{3n}\sqrt{\log n}
\bigr)-\delta_{1n}
\nonumber
\\
&&\qquad \geq\Pr_f\bigl(W_{n,f}^0\leq
c_{n,f}(\alpha+\tau_n)\bigr)-p_{\bar{\epsilon
}_n}\bigl(|G_{n,f}|\bigr)-
\delta_{1n} \label{eqprop1der6}
\\
&&\qquad \geq1-\alpha-\tau_n-p_{\bar{\epsilon}_n}\bigl(|G_{n,f}|\bigr)-
\delta_{1n}, \label{eqprop1der7}
\end{eqnarray}
where (\ref{eqprop1der6}) follows from the definition of the L\'{e}vy
concentration function $p_{\bar{\epsilon}_n}(|G_{n,f}|)$ given that
$\bar{\epsilon}_n=\epsilon_{1n}+\epsilon_{2n}+\epsilon
_{3n}(c_{n,f}(\alpha)+u_n\sqrt{\log n})$, and (\ref{eqprop1der7})
follows since $c_{n,f}(\cdot)$ is the quantile function of $W_{n,f}^0$.
Combining these inequalities leads to (\ref{eqresulthonestnonasymptotic}).

To prove (\ref{eqresultweak}) and (\ref{eqresultstrong}), note that
$\delta_n\leq Cn^{-c}$ and $\tau_n\leq Cn^{-c}$ by Conditions~\ref{con1} and~\ref{con3}--\ref{con6}. Further, it follows from Lemma
\ref{lemcriticalvaluebound} that $c_{n,f}(\alpha)\leq\Ep[\Vert
G_{n,f}\Vert_{\mathcal{V}_n}]+(2|\log\alpha|)^{1/2}\leq C\sqrt
{\log n}$, and so $\epsilon_{3n}u_n\sqrt{\log n} \leq C_1n^{-c_1}$
implies that $\bar{\epsilon}_{n,f}\leq Cn^{-c}$. Therefore, (\ref
{eqresultweak}) and (\ref{eqresultstrong}) follow from (\ref
{eqresulthonestnonasymptotic}) and Condition~\ref{con2}.
\end{pf*}

\begin{pf*}{Proof of Corollary~\ref{corundersmoothing}}
The proof is similar to that of Theorem~\ref{prop1}. The details are
provided in the supplemental material \cite{CCK2014}.
\end{pf*}

\begin{pf*}{Proof of Theorem~\ref{propverifyingconditions1-3}}
In this proof, $c, C>0$ are constants that depend only on $c_2,C_2$,
but their values can change at each appearance.

Fix any $f\in\mathcal{F}$. Let $G_{n,f}=\{G_{n,f}(v)\dvtx v\in\mathcal
{V}_n\}$ be a tight Gaussian random element in $\ell^{\infty
}(\mathcal{V}_{n})$ with mean zero and the same covariance function as
that of $Z_{n,f}$. Since $b_n^2\sigma_n^4K_n^4/n\leq C_2n^{-c_2}$, it
follows from Theorem~\ref{thmcouplinginequality1} that we can
construct a random variable $W_{n,f}^0$ such that $W_{n,f}^0\stackrel
{d}{=}\|G_{n,f}\|_{\mathcal{V}_n}$, and (\ref{eqcon1}) holds with
some $\epsilon_{1n}$ and $\delta_{1n}$ bounded from above by
$Cn^{-c}$. In addition, inequality $\Ep[\|G_{n,f}\|_{\mathcal
{V}_n}]\leq C\sqrt{\log n}$ follows from Corollary 2.2.8 in \cite
{VW96}. Condition~\ref{con1} follows. Given Condition~\ref{con1},
Condition~\ref{con2}(b) follows from Corollary~\ref{cor1}, and
Condition~\ref{con2}(a) follows from Condition~\ref{con2}(b).

Consider Condition~\ref{con4}. There exists $n_0$ such that
$C_2n_0^{-c_2}\leq1$. It suffices to verify the condition only for
$n\geq n_0$. Note that
%
%
\begin{equation}
\label{eqvariancebound1} \biggl\llvert \frac{\hat{\sigma}_n(x,l)}{\sigma_{n,f}(x,l)}-1\biggr\rrvert \leq \biggl\llvert
\frac{\hat{\sigma}^2_n(x,l)}{\sigma^2_{n,f}(x,l)}-1\biggr\rrvert.
\end{equation}
Define $\mathcal{K}^2_{n,f}:=\{g^2\dvtx g\in\mathcal{K}_{n,f}\}$. Given
the definition of $\hat{\sigma}_n(x,l)$, the right-hand side of (\ref
{eqvariancebound1}) is bounded by
%
%
\begin{equation}
\label{eqRHSbound1} \sup_{g\in\mathcal{K}^2_{n,f}}\bigl\llvert \mathbb{E}_n
\bigl[g(X_i)\bigr]-\Ep \bigl[g(X_1)\bigr]\bigr\rrvert +
\sup_{g\in\mathcal{K}_{n,f}}\bigl\llvert \mathbb {E}_n
\bigl[g(X_i)\bigr]^2-\Ep\bigl[g(X_1)
\bigr]^2\bigr\rrvert.
\end{equation}
It follows from Lemma~\ref{lemvcproperties} that $\mathcal
{K}_{n,f}^2$ is VC$(b_n^2,2a_n,v_n)$ type class. Moreover, for all
$g\in\mathcal{K}^2_{n,f}$,
\[
\Ep\bigl[g(X_i)^2\bigr]\leq b_n^2
\Ep\bigl[g(X_i)\bigr]\leq b_n^2
\sigma_n^2.
\]
Therefore, Talagrand's inequality (Theorem~\ref
{thmtalagrandinequality}) with $t=\log n$, which can be applied because\vspace*{1pt}
$b_n^2K_n/(n\sigma_n^2)\leq b_n^2\sigma_n^4K_n^4/n \leq
C_2n^{-c_2}\leq1$ and $b_n^2\log n/\break(n\sigma_n^2)\leq
b_n^2K_n/ (n\sigma_n^2)\leq1$ (recall that $\sigma_n\geq1$ and
$K_n\geq1$), gives
%
%
\begin{equation}
\label{eqfirsttermbound1} \Pr \biggl(\sup_{g\in\mathcal{K}^2_{n,f}}\bigl\llvert \mathbb
{E}_n\bigl[g(X_i)\bigr]-\Ep\bigl[g(X_1)
\bigr]\bigr\rrvert >\frac{1}{2}\sqrt{\frac
{b_n^2\sigma_n^2K_n}{n}} \biggr)
\leq\frac{1}{n}.
\end{equation}
In addition,
\[
\sup_{g\in\mathcal{K}_{n,f}}\bigl\llvert \mathbb{E}_n
\bigl[g(X_i)\bigr]^2-\Ep \bigl[g(X_1)
\bigr]^2\bigr\rrvert \leq2b_n\sup_{g\in\mathcal{K}_{n,f}}
\bigl\llvert \mathbb{E}_n\bigl[g(X_i)\bigr]-\Ep
\bigl[g(X_1)\bigr]\bigr\rrvert,
\]
so that another application of Talagrand's inequality yields
%
%
\begin{equation}
\label{eqsecondtermbound1} \Pr \biggl(\sup_{g\in\mathcal{K}_{n,f}}\bigl\llvert \mathbb
{E}_n\bigl[g(X_i)\bigr]^2-\Ep
\bigl[g(X_1)\bigr]^2\bigr\rrvert >\frac{1}{2}\sqrt
{\frac
{b_n^2\sigma_n^2K_n}{n}} \biggr)\leq\frac{1}{n}.
\end{equation}
Given that $b_n^2\sigma_n^2K_n/n\leq b_n^2\sigma_n^4K_n^4/n\leq
C_2n^{-c_2}$, combining (\ref{eqvariancebound1})--(\ref
{eqsecondtermbound1}) gives Condition~\ref{con4} with $\epsilon
_{3n}:=(b_n^2\sigma_n^2K_n/n)^{1/2}$ and $\delta_{3n}:=2/n$.

Finally, we verify Condition~\ref{con3}. There exists $n_1$ such that
$\epsilon_{3n_1}\leq1/2$. It suffices to verify the condition only
for $n\geq n_1$, so that $\epsilon_{3n}\leq1/2$. Define
\[
\widetilde{\mathbb{G}}_n(x,l)=\widetilde{\mathbb{G}}_n
\bigl(X_1^n,\xi _1^n\bigr)
(x,l):=\frac{1}{\sqrt{n}}\sum_{i=1}^n
\xi_i\frac{K_l(X_i,
x)-\hat{f}_n(x,l)}{\sigma_n(x,l)}
\]
and
\[
\Delta\mathbb{G}_n(x,l)=\hat{\mathbb{G}}_n(x,l)-\widetilde{
\mathbb{G}}_n(x,l).
\]
In addition, define
\begin{eqnarray*}
\widehat{W}_n\bigl(x_1^n\bigr)&:=&\sup
_{(x,l)\in\mathcal{X}\times\mathcal
{L}_n}\hat{\mathbb{G}}_n\bigl(x_1^n,
\xi_1^n\bigr) (x,l),
\\
\widetilde{W}_n\bigl(x_1^n\bigr)&:=&\sup
_{(x,l)\in\mathcal{X}\times\mathcal
{L}_n}\widetilde{\mathbb{G}}_n\bigl(x_1^n,
\xi_1^n\bigr) (x,l).
\end{eqnarray*}
Consider the set $S_{n,1}$ of values $X_1^n$ such that $|\hat{\sigma
}_n(x,l)/\sigma_{n,f}(x,l)-1|\leq\epsilon_{3n}$ for all $(x,l)\in
\mathcal{X}\times\mathcal{L}_n$ whenever $X_1^n\in S_{n,1}$. The
previous calculations show that $\Pr_f(X_1^n \in S_{n,1})\geq1-\delta
_{3n}=1-2/n$. Pick and fix any $x_1^n\in S_{n,1}$. Then
\[
\Delta\mathbb{G}_n\bigl(x_1^n,
\xi_1^n\bigr) (x,l)=\frac{1}{\sqrt{n}}\sum
_{i=1}^n\xi_i\frac{K_l(x_i,x)-\hat{f}_n(x,l)}{\sigma_n(x,l)} \biggl(
\frac{\sigma_n(x,l)}{\hat{\sigma}_n(x,l)}-1 \biggr)
\]
is a Gaussian process with mean zero and
\[
\var \bigl(\Delta\mathbb{G}_n\bigl(x_1^n,
\xi_1^n\bigr) (x,l) \bigr)=\frac
{\hat{\sigma}_n^2(x,l)}{\sigma_n^2(x,l)} \biggl(
\frac{\sigma
_n(x,l)}{\hat{\sigma}_n(x,l)}-1 \biggr)^2\leq\epsilon_{3n}^2.
\]
Further, the function class
\[
\widetilde{\mathcal{K}}_{n,f}:= \biggl\{\frac{K_l(\cdot, x)}{\sigma
_n(x,l)} \biggl(
\frac{\sigma_n(x,l)}{\hat{\sigma}_n(x,l)}-1 \biggr)\dvtx (x,l)\in\mathcal{X}\times\mathcal{L}_n
\biggr\}
\]
is contained in the function class
\[
\biggl\{\frac{aK_l(\cdot,x)}{\sigma_n(x,l)}\dvtx (x,l,a)\in\mathcal {X}\times\mathcal{L}_n
\times[-1,1] \biggr\},
\]
and hence is VC$(b_n,4a_n,1+v_n)$ type class by Lemma~\ref{lemvcproperties}.
In addition,
\begin{eqnarray*}
&& \Ep \bigl[ \bigl(\Delta\mathbb{G}_n\bigl(x_1^n,
\xi_1^n\bigr) \bigl(x',l'
\bigr)-\Delta \mathbb{G}_n\bigl(x_1^n,
\xi_1^n\bigr) \bigl(x'',l''
\bigr) \bigr)^2 \bigr]
\\
&&\qquad \leq \mathbb{E}_n \biggl[ \biggl(\frac{K_l(x_i,x')}{\sigma_n(x',l')} \biggl(
\frac{\sigma_n(x',l')}{\hat{\sigma}_n(x',l')}-1 \biggr)-\frac
{K_l(x_i,x'')}{\sigma_n(x'',l'')} \biggl(\frac{\sigma
_n(x'',l'')}{\hat{\sigma}_n(x'',l'')}-1
\biggr) \biggr)^2 \biggr],
\end{eqnarray*}
for all $x',x''\in\mathcal{X}$ and $l',l''\in\mathcal{L}_n$, so
that covering numbers for the index set \mbox{$\mathcal{X} \times\mathcal
{L}_{n}$} with respect to the intrinsic (standard deviation) semimetric
induced from the Gaussian process $\Delta\mathbb{G}_n(x_1^n,\xi
_1^n)$ are bounded by uniform covering numbers for the function class
$\widetilde{\mathcal{K}}_{n,f}$.
Therefore, an application of Corollary 2.2.8 in \cite{VW96} gives
\begin{eqnarray*}
\Ep \Bigl[ \sup_{(x,l)\in\mathcal{X}\times\mathcal{L}_n}\bigl|\Delta \mathbb{G}_n
\bigl(x_1^n,\xi_1^n\bigr) (x,l)\bigr|
\Bigr]&\leq& C\epsilon_{3n}\sqrt {(1+v_n)\log \biggl(
\frac{4a_nb_n}{\epsilon_{3n}} \biggr)}\leq C n^{-c}.
\end{eqnarray*}
Here the second inequality follows from the definition of $\epsilon
_{3n}$ above and the following inequalities:
\begin{eqnarray*}
\sqrt{(1+v_n)\log \biggl(\frac{4a_n b_n}{\epsilon_{3n}} \biggr)} &\leq&
\sqrt{(1+v_n) \biggl(\log \biggl(\frac{4a_n b_n}{\sigma
_n} \biggr)+\log
\biggl(\frac{\sigma_n}{\epsilon_{3n}} \biggr) \biggr)}
\\
&\leq& C\sqrt{K_n} \biggl(1+\sqrt{\log \biggl(\frac{\sigma
_n}{\epsilon_{3n}}
\biggr)} \biggr)
\\
&\leq& C\sqrt{K_n} \biggl(1+\sqrt{\log \biggl(
\frac{n}{b_n^2
K_n} \biggr)} \biggr)
\\
&\leq& C\sqrt{K_n} (1+\sqrt{\log n} )\leq CK_n,
\end{eqnarray*}
where in the last line we used $b_n\geq\sigma_n\geq1$, and $K_n\geq
v_n\log n\geq\log n$.
Combining this bound with the Borell--Sudakov--Tsirel'son inequality,
and using the inequality
\[
\bigl|\widehat{W}_n\bigl(x_1^n\bigr)-
\widetilde{W}_n\bigl(x_1^n\bigr)\bigr|\leq\sup
_{(x,l)\in\mathcal
{X}\times\mathcal{L}_n}\bigl|\Delta\mathbb{G}_n\bigl(x_1^n,
\xi_1^n\bigr) (x,l)\bigr|,
\]
we see that there exists $\lambda_{1n} \leq Cn^{-c}$ such that
%
%
\begin{equation}
\label{eqlambdaexistence1} \Pr\bigl(\bigl|\widehat{W}_n\bigl(x_1^n
\bigr)-\widetilde{W}_n\bigl(x_1^n\bigr)\bigr|\geq
\lambda_{1n}\bigr)\leq Cn^{-c},
\end{equation}
whenever $x_1^n\in S_{n,1}$.
Further, since $b_n^2\sigma_n^4K_n^4/n\leq C_2n^{-c_2}$ and $b_n\geq
\sigma_n\geq1$, Theorem~\ref{thmcouplinginequality2} shows that
there exist $\lambda_{2n} \leq Cn^{-c}$ and a measurable set $S_{n,2}$
of values $X_1^n$ such that $\Pr_{f} (X_{1}^{n} \in S_{n,2})\geq
1-3/n$, and for every $x_1^n\in S_{n,2}$ one can construct a random
variable $W^0$ such that $W^0\stackrel{d}{=}\Vert G_{n,f}\Vert
_{\mathcal{V}_n}$ and
%
%
\begin{equation}
\label{eqlambdaexistence2} \Pr\bigl(\bigl|\widetilde{W}_n\bigl(x_1^n
\bigr)-W^0\bigr|\geq\lambda_{2n}\bigr)\leq Cn^{-c}.
\end{equation}
Here $W^{0}$ may depend on $x_{1}^{n}$, but $c,C$ can be chosen in such
a way that they depend only on $c_{2},C_{2}$ (as noted in the beginning).

Pick and fix any $x_{1}^{n} \in S_{n,0}:=S_{n,1}\cap S_{n,2}$, and
construct a suitable $W^{0} \stackrel{d}{=} \Vert G_{n,f}\Vert
_{\mathcal{V}_n}$ for which (\ref{eqlambdaexistence2}) holds.
Then by (\ref{eqlambdaexistence1}), we have
%
%
\begin{equation}
\label{eqWnfandW0close} \Pr\bigl(\bigl|\widehat{W}_n\bigl(x_1^n
\bigr)-W^0\bigr|\geq\lambda_{n} \bigr)\leq Cn^{-c},
\end{equation}
where $\lambda_{n}:= \lambda_{1n} + \lambda_{2n}$.
Denote by $\hat{c}_n(\alpha,x_1^n)$ the $(1-\alpha)$-quantile of
$\widehat{W}_n(x_1^n)$. Then we have
%
\begin{eqnarray*}
\Pr\bigl( \| G_{n,f} \|_{\mathcal{V}_{n}} \leq\hat{c}_{n}\bigl(
\alpha,x_1^n\bigr) + \lambda_n \bigr)&=&\Pr
\bigl(W^0 \leq\hat{c}_{n}\bigl(\alpha,x_1^n
\bigr) + \lambda_n \bigr)
\\
&\geq&\Pr\bigl(\widehat{W}_n\bigl(x_1^n
\bigr)\leq\hat{c}_{n}\bigl(\alpha,x_1^n\bigr)
\bigr) - Cn^{-c}
\\
&\geq&1-\alpha- Cn^{-c},
\end{eqnarray*}
by which we have $\hat{c}_{n}(\alpha,x_{1}^{n}) \geq c_{n,f}(\alpha+
Cn^{-c}) - \lambda_n$. Since $x_{1}^{n} \in S_{n,0}$ is arbitrary and
$\hat{c}_n(\alpha)=\hat{c}_n(\alpha,X_1^n)$, we see that whenever
$X_{1}^{n} \in S_{n,0}$,
$\hat{c}_{n}(\alpha) \geq c_{n,f}(\alpha+ Cn^{-c}) - \lambda_n$.
Part (a) of Condition~\ref{con3} follows from the fact that $\Pr
_{f}(X_{1}^{n} \in S_{n,0}) \geq1-5/n$ and $\lambda_n \leq Cn^{-c}$.
Part (b) follows similarly.
\end{pf*}
\end{appendix}

\section*{Acknowledgments}
The authors would like to thank Emre Barut, Enno Mammen and Richard
Nickl for very helpful discussions. We also thank the editors and
anonymous referees for their very helpful reviews that helped improve
the paper.

\begin{supplement}[id=suppA]
\stitle{Supplement to ``Anti-concentration and honest, adaptive confidence bands''}
\slink[doi]{10.1214/14-AOS1235SUPP} 
\sdatatype{.pdf}
\sfilename{aos1235\_supp.pdf}
\sdescription{This supplemental file contains additional proofs
omitted in the main text, some results regarding nonwavelet projection
kernel estimators, and a small-scale simulation study.}
\end{supplement}


\printaddresses

\begin{thebibliography}{45}

\bibitem{BR73}
\begin{barticle}[mr]
\bauthor{\bsnm{Bickel},~\bfnm{P.~J.}\binits{P.~J.}} \AND
\bauthor{\bsnm{Rosenblatt},~\bfnm{M.}\binits{M.}}
(\byear{1973}).
\btitle{On some global measures of the deviations of density function estimates}.
\bjournal{Ann. Statist.}
\bvolume{1}
\bpages{1071--1095}.
\bid{issn={0090-5364}, mr={0348906}}%
\end{barticle}%
\bptok{imsref}%
\endbibitem

\bibitem{Birge2001}
\begin{bincollection}[mr]
\bauthor{\bsnm{Birg{\'e}},~\bfnm{Lucien}\binits{L.}}
(\byear{2001}).
\btitle{An alternative point of view on {L}epski's method}.
In \bbooktitle{State of the Art in Probability and Statistics ({L}eiden, 1999)}.
\bseries{Institute of Mathematical Statistics Lecture Notes---Monograph Series}
\bvolume{36}
\bpages{113--133}.
\bpublisher{IMS},
\blocation{Beachwood, OH}.
\bid{doi={10.1214/lnms/1215090065}, mr={1836557}}
\end{bincollection}
\bptok{imsref}%
\endbibitem

\bibitem{BDHM07}
\begin{barticle}[mr]
\bauthor{\bsnm{Bissantz},~\bfnm{Nicolai}\binits{N.}},
\bauthor{\bsnm{D{\"u}mbgen},~\bfnm{Lutz}\binits{L.}},
\bauthor{\bsnm{Holzmann},~\bfnm{Hajo}\binits{H.}} \AND
\bauthor{\bsnm{Munk},~\bfnm{Axel}\binits{A.}}
(\byear{2007}).
\btitle{Non-parametric confidence bands in deconvolution density estimation}.
\bjournal{J. R. Stat. Soc. Ser. B Stat. Methodol.}
\bvolume{69}
\bpages{483--506}.
\bid{doi={10.1111/j.1467-9868.2007.599.x}, issn={1369-7412}, mr={2323764}}
\end{barticle}
\bptok{imsref}%
\endbibitem

\bibitem{B11a}
\begin{barticle}[mr]
\bauthor{\bsnm{Bull},~\bfnm{Adam~D.}\binits{A.~D.}}
(\byear{2012}).
\btitle{Honest adaptive confidence bands and self-similar functions}.
\bjournal{Electron. J.~Stat.}
\bvolume{6}
\bpages{1490--1516}.
\bid{doi={10.1214/12-EJS720}, issn={1935-7524}, mr={2988456}}
\end{barticle}
\bptok{imsref}%
\endbibitem

\bibitem{B11b}
\begin{barticle}[mr]
\bauthor{\bsnm{Bull},~\bfnm{Adam~D.}\binits{A.~D.}}
(\byear{2013}).
\btitle{A {S}mirnov--{B}ickel--{R}osenblatt theorem for compactly-supported wavelets}.
\bjournal{Constr. Approx.}
\bvolume{37}
\bpages{295--309}.
\bid{doi={10.1007/s00365-013-9181-7}, issn={0176-4276}, mr={3019781}}
\end{barticle}
\bptok{imsref}%
\endbibitem

\bibitem{CaiLow04}
\begin{barticle}[mr]
\bauthor{\bsnm{Cai},~\bfnm{T.~Tony}\binits{T.~T.}} \AND
\bauthor{\bsnm{Low},~\bfnm{Mark~G.}\binits{M.~G.}}
(\byear{2004}).
\btitle{An adaptation theory for nonparametric confidence intervals}.
\bjournal{Ann. Statist.}
\bvolume{32}
\bpages{1805--1840}.
\bid{doi={10.1214/009053604000000049}, issn={0090-5364}, mr={2102494}}
\end{barticle}
\bptok{imsref}%
\endbibitem

\bibitem{CCK2012}
\begin{bmisc}[auto:STB|2014/06/18|12:29:53]
\bauthor{\bsnm{Chernozhukov},~\bfnm{V.}\binits{V.}},
\bauthor{\bsnm{Chetverikov},~\bfnm{D.}\binits{D.}}
\and
\bauthor{\bsnm{Kato},~\bfnm{K.}\binits{K.}}
(\byear{2012}).
\bhowpublished{Gaussian approximation of suprema of empirical processes.
Preprint. Available at \arxivurl{arXiv:1212.6885v2}.}
\end{bmisc}
\bptok{imsref}%
\endbibitem

\bibitem{CCK2012b}
\begin{bmisc}[auto:STB|2014/06/18|12:29:53]
\bauthor{\bsnm{Chernozhukov},~\bfnm{V.}\binits{V.}},
\bauthor{\bsnm{Chetverikov},~\bfnm{D.}\binits{D.}}
\and
\bauthor{\bsnm{Kato},~\bfnm{K.}\binits{K.}}
(\byear{2014}).
\bhowpublished{Comparison and anti-concentration bounds for maxima of Gaussian random vectors.
\textit{Probab. Theory Related Fields}. To appear. Available at \arxivurl{arXiv:1301.4807v3}.}
\end{bmisc}
\bptok{imsref}%
\endbibitem

\bibitem{CCK2014}
\begin{bmisc}[auto:STB|2014/06/18|12:29:53]
\bauthor{\bsnm{Chernozhukov},~\bfnm{V.}\binits{V.}},
\bauthor{\bsnm{Chetverikov},~\bfnm{D.}\binits{D.}}
\and
\bauthor{\bsnm{Kato},~\bfnm{K.}\binits{K.}}
(\byear{2014}).
\bhowpublished{Supplement to ``Anti-concentration and honest, adaptive confidence bands.''
 DOI:\doiurl{10.1214/14-AOS1235SUPP}.}
\end{bmisc}
\bptok{imsref}%
\endbibitem

\bibitem{CLR12}
\begin{barticle}[mr]
\bauthor{\bsnm{Chernozhukov},~\bfnm{Victor}\binits{V.}},
\bauthor{\bsnm{Lee},~\bfnm{Sokbae}\binits{S.}} \AND
\bauthor{\bsnm{Rosen},~\bfnm{Adam~M.}\binits{A.~M.}}
(\byear{2013}).
\btitle{Intersection bounds: Estimation and inference}.
\bjournal{Econometrica}
\bvolume{81}
\bpages{667--737}.
\bid{doi={10.3982/ECTA8718}, issn={0012-9682}, mr={3043345}}
\end{barticle}
\bptok{imsref}%
\endbibitem

\bibitem{Chetverikov2012}
\begin{bmisc}[auto:STB|2014/06/18|12:29:53]
\bauthor{\bsnm{Chetverikov},~\bfnm{D.}\binits{D.}}
(\byear{2012}).
\bhowpublished{Testing regression monotonicity in econometric models.
Available at \arxivurl{arXiv:1212.6757}.}
\end{bmisc}
\bptok{imsref}%
\endbibitem

\bibitem{CV03}
\begin{barticle}[mr]
\bauthor{\bsnm{Claeskens},~\bfnm{Gerda}\binits{G.}} \AND
\bauthor{\bsnm{Van Keilegom},~\bfnm{Ingrid}\binits{I.}}
(\byear{2003}).
\btitle{Bootstrap confidence bands for regression curves and their derivatives}.
\bjournal{Ann. Statist.}
\bvolume{31}
\bpages{1852--1884}.
\bid{doi={10.1214/aos/1074290329}, issn={0090-5364}, mr={2036392}}
\end{barticle}
\bptok{imsref}%
\endbibitem

\bibitem{Daubechies92}
\begin{bbook}[auto:STB|2014/06/18|12:29:53]
\bauthor{\bsnm{Daubechies},~\bfnm{I.}\binits{I.}}
(\byear{1992}).
\btitle{Ten Lectures on Wavelets}.
\bseries{CBMS-NSF Regional Conference Series in Applied Mathematics}
\bvolume{61}.
\bpublisher{SIAM},
\blocation{Philadelphia, PA}.
\bid{mr={1162107}}
\end{bbook}
\bptok{imsref}%
\endbibitem

\bibitem{D99}
\begin{bbook}[mr]
\bauthor{\bsnm{Dudley},~\bfnm{R.~M.}\binits{R.~M.}}
(\byear{1999}).
\btitle{Uniform Central Limit Theorems}.
\bpublisher{Cambridge Univ. Press},
\blocation{Cambridge}.
\bid{doi={10.1017/CBO9780511665622}, mr={1720712}}
\end{bbook}
\bptok{imsref}%
\endbibitem

\bibitem{GG01}
\begin{barticle}[auto]
\bauthor{\bsnm{Gin{\'e}},~\bfnm{Evarist}\binits{E.}} \AND
\bauthor{\bsnm{Guillou},~\bfnm{Armelle}\binits{A.}}
(\byear{2001}).
\btitle{A law of the iterated logarithm for kernel density estimators in the presence of censoring}.
\bjournal{Ann. Inst. Henri Poincar\'e Probab. Stat.}
\bvolume{37}
\bpages{503--522}.
\end{barticle}
\bptok{imsref}%
\endbibitem

\bibitem{GG02}
\begin{barticle}[mr]
\bauthor{\bsnm{Gin{\'e}},~\bfnm{Evarist}\binits{E.}} \AND
\bauthor{\bsnm{Guillou},~\bfnm{Armelle}\binits{A.}}
(\byear{2002}).
\btitle{Rates of strong uniform consistency for multivariate kernel density estimators}.
\bjournal{Ann. Inst. Henri Poincar\'e Probab. Stat.}
\bvolume{38}
\bpages{907--921}.
\bid{doi={10.1016/S0246-0203(02)01128-7}, issn={0246-0203}, mr={1955344}}
\end{barticle}
\bptok{imsref}%
\endbibitem

\bibitem{GGM11}
\begin{barticle}[mr]
\bauthor{\bsnm{Gin{\'e}},~\bfnm{E.}\binits{E.}},
\bauthor{\bsnm{G{\"u}nt{\"u}rk},~\bfnm{C.~S.}\binits{C.~S.}} \AND
\bauthor{\bsnm{Madych},~\bfnm{W.~R.}\binits{W.~R.}}
(\byear{2011}).
\btitle{On the periodized square of {$L^2$} cardinal splines}.
\bjournal{Exp. Math.}
\bvolume{20}
\bpages{177--188}.
\bid{doi={10.1080/10586458.2011.564543}, issn={1058-6458}, mr={2821389}}
\end{barticle}
\bptok{imsref}%
\endbibitem

\bibitem{GKS04}
\begin{barticle}[mr]
\bauthor{\bsnm{Gin{\'e}},~\bfnm{Evarist}\binits{E.}},
\bauthor{\bsnm{Koltchinskii},~\bfnm{Vladimir}\binits{V.}} \AND
\bauthor{\bsnm{Sakhanenko},~\bfnm{Lyudmila}\binits{L.}}
(\byear{2004}).
\btitle{Kernel density estimators: Convergence in distribution for weighted sup-norms}.
\bjournal{Probab. Theory Related Fields}
\bvolume{130}
\bpages{167--198}.
\bid{doi={10.1007/s00440-004-0339-x}, issn={0178-8051}, mr={2093761}}
\end{barticle}
\bptok{imsref}%
\endbibitem

\bibitem{GN09b}
\begin{barticle}[mr]
\bauthor{\bsnm{Gin{\'e}},~\bfnm{Evarist}\binits{E.}} \AND
\bauthor{\bsnm{Nickl},~\bfnm{Richard}\binits{R.}}
(\byear{2009}).
\btitle{An exponential inequality for the distribution function of the kernel density estimator, with applications to adaptive estimation}.
\bjournal{Probab. Theory Related Fields}
\bvolume{143}
\bpages{569--596}.
\bid{doi={10.1007/s00440-008-0137-y}, issn={0178-8051}, mr={2475673}}
\bptnote{check year}
\end{barticle}
\bptok{imsref}%
\endbibitem

\bibitem{GN10}
\begin{barticle}[mr]
\bauthor{\bsnm{Gin{\'e}},~\bfnm{Evarist}\binits{E.}} \AND
\bauthor{\bsnm{Nickl},~\bfnm{Richard}\binits{R.}}
(\byear{2010}).
\btitle{Confidence bands in density estimation}.
\bjournal{Ann. Statist.}
\bvolume{38}
\bpages{1122--1170}.
\bid{doi={10.1214/09-AOS738}, issn={0090-5364}, mr={2604707}}
\end{barticle}
\bptok{imsref}%
\endbibitem

\bibitem{GN10b}
\begin{barticle}[mr]
\bauthor{\bsnm{Gin{\'e}},~\bfnm{Evarist}\binits{E.}} \AND
\bauthor{\bsnm{Nickl},~\bfnm{Richard}\binits{R.}}
(\byear{2010}).
\btitle{Adaptive estimation of a distribution function and its density in sup-norm loss by wavelet and spline projections}.
\bjournal{Bernoulli}
\bvolume{16}
\bpages{1137--1163}.
\bid{doi={10.3150/09-BEJ239}, issn={1350-7265}, mr={2759172}}
\end{barticle}
\bptok{imsref}%
\endbibitem

\bibitem{GZ84}
\begin{barticle}[mr]
\bauthor{\bsnm{Gin{\'e}},~\bfnm{Evarist}\binits{E.}} \AND
\bauthor{\bsnm{Zinn},~\bfnm{Joel}\binits{J.}}
(\byear{1984}).
\btitle{Some limit theorems for empirical processes}.
\bjournal{Ann. Probab.}
\bvolume{12}
\bpages{929--998}.
\bid{issn={0091-1798}, mr={0757767}}
\bptnote{check related}%
\end{barticle}
\bptok{imsref}%
\endbibitem

\bibitem{GZ90}
\begin{barticle}[mr]
\bauthor{\bsnm{Gin{\'e}},~\bfnm{Evarist}\binits{E.}} \AND
\bauthor{\bsnm{Zinn},~\bfnm{Joel}\binits{J.}}
(\byear{1990}).
\btitle{Bootstrapping general empirical measures}.
\bjournal{Ann. Probab.}
\bvolume{18}
\bpages{851--869}.
\bid{issn={0091-1798}, mr={1055437}}
\end{barticle}
\bptok{imsref}%
\endbibitem

\bibitem{H91}
\begin{barticle}[mr]
\bauthor{\bsnm{Hall},~\bfnm{Peter}\binits{P.}}
(\byear{1991}).
\btitle{On convergence rates of suprema}.
\bjournal{Probab. Theory Related Fields}
\bvolume{89}
\bpages{447--455}.
\bid{doi={10.1007/BF01199788}, issn={0178-8051}, mr={1118558}}
\end{barticle}
\bptok{imsref}%
\endbibitem

\bibitem{HH13}
\begin{barticle}[mr]
\bauthor{\bsnm{Hall},~\bfnm{Peter}\binits{P.}} \AND
\bauthor{\bsnm{Horowitz},~\bfnm{Joel}\binits{J.}}
(\byear{2013}).
\btitle{A simple bootstrap method for constructing nonparametric confidence bands for functions}.
\bjournal{Ann. Statist.}
\bvolume{41}
\bpages{1892--1921}.
\bid{doi={10.1214/13-AOS1137}, issn={0090-5364}, mr={3127852}}
\end{barticle}
\bptok{imsref}%
\endbibitem

\bibitem{Hardle}
\begin{bbook}[mr]
\bauthor{\bsnm{H{\"a}rdle},~\bfnm{Wolfgang}\binits{W.}},
\bauthor{\bsnm{Kerkyacharian},~\bfnm{Gerard}\binits{G.}},
\bauthor{\bsnm{Picard},~\bfnm{Dominique}\binits{D.}} \AND
\bauthor{\bsnm{Tsybakov},~\bfnm{Alexander}\binits{A.}}
(\byear{1998}).
\btitle{Wavelets, Approximation, and Statistical Applications}.
\bpublisher{Springer},
\blocation{New York}.
\bid{doi={10.1007/978-1-4612-2222-4}, mr={1618204}}
\end{bbook}
\bptok{imsref}%
\endbibitem

\bibitem{HN11}
\begin{barticle}[mr]
\bauthor{\bsnm{Hoffmann},~\bfnm{Marc}\binits{M.}} \AND
\bauthor{\bsnm{Nickl},~\bfnm{Richard}\binits{R.}}
(\byear{2011}).
\btitle{On adaptive inference and confidence bands}.
\bjournal{Ann. Statist.}
\bvolume{39}
\bpages{2383--2409}.
\bid{doi={10.1214/11-AOS903}, issn={0090-5364}, mr={2906872}}
\end{barticle}
\bptok{imsref}%
\endbibitem

\bibitem{KMT75}
\begin{barticle}[mr]
\bauthor{\bsnm{Koml{\'o}s},~\bfnm{J.}\binits{J.}},
\bauthor{\bsnm{Major},~\bfnm{P.}\binits{P.}} \AND
\bauthor{\bsnm{Tusn{\'a}dy},~\bfnm{G.}\binits{G.}}
(\byear{1975}).
\btitle{An approximation of partial sums of independent {${\rm RV}$}'s and the sample {${\rm DF}$}. {I}}.
\bjournal{Z. Wahrsch. Verw. Gebiete}
\bvolume{32}
\bpages{111--131}.
\bid{mr={0375412}}
\end{barticle}
\bptok{imsref}%
\endbibitem

\bibitem{LT88}
\begin{barticle}[mr]
\bauthor{\bsnm{Ledoux},~\bfnm{M.}\binits{M.}} \AND
\bauthor{\bsnm{Talagrand},~\bfnm{M.}\binits{M.}}
(\byear{1988}).
\btitle{Un crit\`ere sur les petites boules dans le th\'eor\`eme limite central}.
\bjournal{Probab. Theory Related Fields}
\bvolume{77}
\bpages{29--47}.
\bid{doi={10.1007/BF01848129}, issn={0178-8051}, mr={0921817}}
\end{barticle}
\bptok{imsref}%
\endbibitem

\bibitem{Lepski91}
\begin{barticle}[auto]
\bauthor{\bsnm{Lepski{\u\i}},~\bfnm{O.~V.}\binits{O.~V.}}
(\byear{1991}).
\btitle{Asymptotically minimax adaptive estimation. {I}. {U}pper bounds. {O}ptimally adaptive estimates}.
\bjournal{Theory Probab. Appl.}
\bvolume{36}
\bpages{682--697}.
\end{barticle}
\bptok{imsref}%
\endbibitem

\bibitem{L89}
\begin{barticle}[mr]
\bauthor{\bsnm{Li},~\bfnm{Ker-Chau}\binits{K.-C.}}
(\byear{1989}).
\btitle{Honest confidence regions for nonparametric regression}.
\bjournal{Ann. Statist.}
\bvolume{17}
\bpages{1001--1008}.
\bid{doi={10.1214/aos/1176347253}, issn={0090-5364}, mr={1015135}}
\end{barticle}
\bptok{imsref}%
\endbibitem

\bibitem{LN11}
\begin{barticle}[mr]
\bauthor{\bsnm{Lounici},~\bfnm{Karim}\binits{K.}} \AND
\bauthor{\bsnm{Nickl},~\bfnm{Richard}\binits{R.}}
(\byear{2011}).
\btitle{Global uniform risk bounds for wavelet deconvolution estimators}.
\bjournal{Ann. Statist.}
\bvolume{39}
\bpages{201--231}.
\bid{doi={10.1214/10-AOS836}, issn={0090-5364}, mr={2797844}}
\end{barticle}
\bptok{imsref}%
\endbibitem

\bibitem{Low97}
\begin{barticle}[mr]
\bauthor{\bsnm{Low},~\bfnm{Mark~G.}\binits{M.~G.}}
(\byear{1997}).
\btitle{On nonparametric confidence intervals}.
\bjournal{Ann. Statist.}
\bvolume{25}
\bpages{2547--2554}.
\bid{doi={10.1214/aos/1030741084}, issn={0090-5364}, mr={1604412}}
\end{barticle}
\bptok{imsref}%
\endbibitem

\bibitem{M00}
\begin{barticle}[mr]
\bauthor{\bsnm{Massart},~\bfnm{Pascal}\binits{P.}}
(\byear{2000}).
\btitle{About the constants in {T}alagrand's concentration inequalities for empirical processes}.
\bjournal{Ann. Probab.}
\bvolume{28}
\bpages{863--884}.
\bid{doi={10.1214/aop/1019160263}, issn={0091-1798}, mr={1782276}}
\end{barticle}
\bptok{imsref}%
\endbibitem

\bibitem{M07}
\begin{bbook}[mr]
\bauthor{\bsnm{Massart},~\bfnm{Pascal}\binits{P.}}
(\byear{2007}).
\btitle{Concentration Inequalities and Model Selection}.
\bpublisher{Springer},
\blocation{Berlin}.
\bid{mr={2319879}}
\end{bbook}
\bptok{imsref}%
\endbibitem

\bibitem{PT00}
\begin{barticle}[mr]
\bauthor{\bsnm{Picard},~\bfnm{Dominique}\binits{D.}} \AND
\bauthor{\bsnm{Tribouley},~\bfnm{Karine}\binits{K.}}
(\byear{2000}).
\btitle{Adaptive confidence interval for pointwise curve estimation}.
\bjournal{Ann. Statist.}
\bvolume{28}
\bpages{298--335}.
\bid{doi={10.1214/aos/1016120374}, issn={0090-5364}, mr={1762913}}
\end{barticle}
\bptok{imsref}%
\endbibitem

\bibitem{Wellner93}
\begin{barticle}[mr]
\bauthor{\bsnm{Pr{\ae}stgaard},~\bfnm{Jens}\binits{J.}} \AND
\bauthor{\bsnm{Wellner},~\bfnm{Jon~A.}\binits{J.~A.}}
(\byear{1993}).
\btitle{Exchangeably weighted bootstraps of the general empirical process}.
\bjournal{Ann. Probab.}
\bvolume{21}
\bpages{2053--2086}.
\bid{issn={0091-1798}, mr={1245301}}
\end{barticle}
\bptok{imsref}%
\endbibitem

\bibitem{Rio94}
\begin{barticle}[mr]
\bauthor{\bsnm{Rio},~\bfnm{Emmanuel}\binits{E.}}
(\byear{1994}).
\btitle{Local invariance principles and their application to density estimation}.
\bjournal{Probab. Theory Related Fields}
\bvolume{98}
\bpages{21--45}.
\bid{doi={10.1007/BF01311347}, issn={0178-8051}, mr={1254823}}
\end{barticle}
\bptok{imsref}%
\endbibitem

\bibitem{RV09}
\begin{barticle}[mr]
\bauthor{\bsnm{Rudelson},~\bfnm{Mark}\binits{M.}} \AND
\bauthor{\bsnm{Vershynin},~\bfnm{Roman}\binits{R.}}
(\byear{2009}).
\btitle{Smallest singular value of a random rectangular matrix}.
\bjournal{Comm. Pure Appl. Math.}
\bvolume{62}
\bpages{1707--1739}.
\bid{doi={10.1002/cpa.20294}, issn={0010-3640}, mr={2569075}}
\end{barticle}
\bptok{imsref}%
\endbibitem

\bibitem{S50}
\begin{barticle}[mr]
\bauthor{\bsnm{Smirnov},~\bfnm{N.~V.}\binits{N.~V.}}
(\byear{1950}).
\btitle{On the construction of confidence regions for the density of distribution of random variables}.
\bjournal{Doklady Akad. Nauk SSSR (N.S.)}
\bvolume{74}
\bpages{189--191}.
\bid{mr={0037494}}
\end{barticle}
\bptok{imsref}%
\endbibitem

\bibitem{T96}
\begin{barticle}[mr]
\bauthor{\bsnm{Talagrand},~\bfnm{Michel}\binits{M.}}
(\byear{1996}).
\btitle{New concentration inequalities in product spaces}.
\bjournal{Invent. Math.}
\bvolume{126}
\bpages{505--563}.
\bid{doi={10.1007/s002220050108}, issn={0020-9910}, mr={1419006}}
\end{barticle}
\bptok{imsref}%
\endbibitem

\bibitem{Tsybakov09}
\begin{bbook}[mr]
\bauthor{\bsnm{Tsybakov},~\bfnm{Alexandre~B.}\binits{A.~B.}}
(\byear{2009}).
\btitle{Introduction to Nonparametric Estimation}.
\bpublisher{Springer},
\blocation{New York}.
\bid{doi={10.1007/b13794}, mr={2724359}}
\end{bbook}
\bptok{imsref}%
\endbibitem

\bibitem{VW96}
\begin{bbook}[mr]
\bauthor{\bsnm{van~der Vaart},~\bfnm{Aad~W.}\binits{A.~W.}} \AND
\bauthor{\bsnm{Wellner},~\bfnm{Jon~A.}\binits{J.~A.}}
(\byear{1996}).
\btitle{Weak Convergence and Empirical Processes: With Applications to Statistics}.
\bpublisher{Springer},
\blocation{New York}.
\bid{doi={10.1007/978-1-4757-2545-2}, mr={1385671}}
\end{bbook}
\bptok{imsref}%
\endbibitem

\bibitem{Wasserman06}
\begin{bbook}[mr]
\bauthor{\bsnm{Wasserman},~\bfnm{Larry}\binits{L.}}
(\byear{2006}).
\btitle{All of Nonparametric Statistics}.
\bpublisher{Springer},
\blocation{New York}.
\bid{mr={2172729}}
\end{bbook}
\bptok{imsref}%
\endbibitem

\bibitem{X98}
\begin{barticle}[mr]
\bauthor{\bsnm{Xia},~\bfnm{Yingcun}\binits{Y.}}
(\byear{1998}).
\btitle{Bias-corrected confidence bands in nonparametric regression}.
\bjournal{J. R. Stat. Soc. Ser. B Stat. Methodol.}
\bvolume{60}
\bpages{797--811}.
\bid{doi={10.1111/1467-9868.00155}, issn={1369-7412}, mr={1649488}}
\end{barticle}
\bptok{imsref}%
\endbibitem

\end{thebibliography}
\end{document}